\documentstyle[amsmath,amssymb,graphicx]{article}

\begin{document}

\thispagestyle{empty}

\baselineskip15pt

\title{{\bf Discriminants of Symmetric Polynomials
} \vspace{.2cm}}
\author{{\bf N.Perminov}\footnote{ {\small {\it
Department of General Relativity and Gravitation,
Kazan State University, Kremlevskaya str. 18, Kazan 420008, Russia}};
Nikolay.Perminov@ksu.ru} \ and {\bf
Sh.Shakirov}\thanks{{\small {\it ITEP, Moscow, Russia} and {\it MIPT, Dolgoprudny, Russia}};
shakirov@itep.ru} \date{ }}

\maketitle

\vspace{-4.5cm}

\begin{center}
\hfill ITEP/TH-53/09
\end{center}

\vspace{3.4cm}

\centerline{ABSTRACT}

\bigskip

{\footnotesize
A homogeneous polynomial $S(x_1, \ldots, x_n)$ of degree $r$ in $n$ variables posesses a discriminant $D_{n|r}(S)$, which vanishes if and only if the system of equations $\partial S / \partial x_i = 0$ has non-trivial solutions. We give an explicit formula for discriminants of symmetric (under permutations of $x_1, \ldots, x_n$) homogeneous polynomials of degree $r$ in $n \geq r$ variables. This formula is division free and quite effective from the computational point of view: symbolic computer calculations with the help of this formula take seconds even for $n \approx 20$. We work out in detail the cases $r = 2,3,4$ which will be probably important in applications. We also consider the case of completely antisymmetric polynomials. }


\bigskip

\tableofcontents

\section{Introduction}

This paper deals with discriminants of homogeneous polynomials $S(x_1, \ldots, x_n)$ of degree $r$ in $n$ variables. Discriminants are classical objects of study both in algebra and in geometry. The foundations of discriminant theory were laid in the 18th century by Cayley, Sylvester and Bezout \cite{Foundations}. In this period mostly the case of $n = 2$ variables was considered, with some rare exceptions. The theory of discriminants for $n > 2$ and closely related objects, namely resultants and hyperdeterminants, was significantly developed by Gelfand, Kapranov and Zelevinsky in their seminal book \cite{GKZ}. Recently, a suggestion and attempt to study discriminants and resultants from the viewpoint of quantum field and string theory was made by Dolotin and Morozov \cite{NOLINAL}. This direction of research, which is called \emph{non-linear algebra}, poses a lot of new questions and reveals very promising connections between path integrals and algebraic geometry \cite{INTDISC}. This paper can be considered as yet another contribution to the non-linear algebra program.

Actually, discriminants naturally appear in many physical contexts, for example, in classical gravity. In Riemannian geometry, the main object of study is a smooth manifold equipped with a quadratic form

\begin{align*}
 ds^2 = g_{ij} dx_i dx_j
\end{align*}
\smallskip\\
which is called the metric and often required to be non-degenerate. Non-degeneracy of the metric is fully controlled by its determinant $\det(g)$: if the latter vanishes, the form becomes degenerate. As we know, this metric structure stands behind the classical Einstein theory of relativity, rich enough to describe a variety of observable phenomena in the Solar System and beyond. A still richer theory one can obtain by replacing the quadratic form by a form of any higher degree:

\begin{align*}
ds^3 = g_{ijk} dx_i dx_j dx_k
\end{align*}

\begin{align*}
ds^4 = g_{ijkl} dx_i dx_j dx_k dx_l
\end{align*}
\smallskip\\
and so on. The corresponding geometry is known as Finslerian geometry (for a review see \cite{Rund} and references therein). This far-going generalisation of Riemannian geometry allows to describe even more phenomena and has certain applications in cosmology \cite{Cosm}. The generalised metric $g$ has more than two indices, i.e, it is no longer a matrix and does not posess a determinant in the usual sence. To determine its non-degeneracy, one needs an analogue of the determinant for higher degree forms: the discriminant. For more information about applications of discriminants and resultants to Finslerian geometry, see \cite{Perm}.

To give a more precise definition of discriminant, let us recall the closely related and somewhat more general notion of \emph{resultant}. Consider a system of $n$ equations

\[
\left\{
\begin{array}{ccc}
f_1(x_1, \ldots, x_n) = 0\\
\\
f_2(x_1, \ldots, x_n) = 0\\
\\
\ldots\\
\\
f_n(x_1, \ldots, x_n) = 0\\
\\
\end{array}
\right.
\]
\smallskip\\
where $f_1, \ldots, f_n$ are homogeneous polynomials of degrees $r_1, \ldots, r_n$. This system of $n$ polynomial equations in $n$ variables is over-defined: in general position it does not have non-zero solutions $(x_1, \ldots, x_n)$ at all. To have a non-zero solution, its coefficients should satisfy a single algebraic constraint

\begin{align}
 R \big\{ f_1, f_2, \ldots, f_n \big\} = 0
\end{align}
\smallskip\\
where $R$ is an irreducible polynomial function called resultant, depending on coefficients of the non-linear system under consideration. The existence of such a function, which defines a solvability condition for non-linear systems, can be proved \cite{GKZ}. It can be also proved that resultant is unique up to overall constant factor (which is not important for most applications). Another classical theorem is that resultant is a homogeneous polynomial in coefficients of the system of degree

\begin{align}
\deg R \big\{ f_1, f_2, \ldots, f_n \big\} = r_1 r_2 \ldots r_n \left( \dfrac{1}{r_1} + \dfrac{1}{r_2} + \ldots + \dfrac{1}{r_n} \right)
\label{DegRes}
\end{align}
\smallskip\\
Discriminant is a particular case of resultant, associated with systems of equations of the form

\[
\left\{
\begin{array}{ccc}
\dfrac{\partial S}{\partial x_1}(x_1, \ldots, x_n) = 0\\
\\
\dfrac{\partial S}{\partial x_2}(x_1, \ldots, x_n) = 0\\
\\
\ldots\\
\\
\dfrac{\partial S}{\partial x_n}(x_1, \ldots, x_n) = 0\\
\\
\end{array}
\right.
\]
\smallskip\\
where $S(x_1, \ldots, x_n)$ is a homogeneous polynomial of degree $r$ in $n$ variables. Non-zero solutions $(x_1, \ldots, x_n)$ of this system are called \emph{critical points}, its resultant is called the \emph{discriminant} of $S$ and denoted as $D_{n|r}(S)$:

\begin{align}
D_{n|r}(S) = R\left\{ \dfrac{\partial S}{\partial x_1}, \ldots, \dfrac{\partial S}{\partial x_n} \right\}
\label{Discrim}
\end{align}
\smallskip\\
As follows from eq. (\ref{DegRes}), discriminant has degree

\begin{align}
\deg D_{n|r}(S) = n(r-1)^{n-1}
\label{DegDisc}
\end{align}
\smallskip\\
in coefficients of $S$. In this paper, we fix the overall normalisation of $D_{n|r}(S)$ by requiring that

\begin{align}
D_{n|r}\big(x_1^r + \ldots + x_n^r\big) = 1
\label{Normalisation}
\end{align}
\smallskip\\
In the case of quadratic forms, i.e, for $r = 2$, discriminant is nothing but the determinant $D_{n|2}(S) = \det(S)$. For $r > 2$ the discriminant is more complicated. It is not usually straightforward to calculate a discriminant in a practical situation. Analogues of techniques known for determinants (say, expansion in minors or the "log det = trace log" formula) are either not available or less efficient for discriminants.

In view of this, it becomes especially important to find formulas which have practical significance, i.e, allow fast and efficient calculation of discriminants, at least in some particular cases. In this paper, we study one such particular case -- of polynomials $S(x_1, \ldots, x_n)$ which are symmetric under permutations of $x_1, \ldots, x_n$. It appears that restriction to symmetric polynomials simplifies the discriminant and allows to calculate it for arbitrarily large number $n$ of variables, see eq. (\ref{DegreeTwo}) and especially eq. (\ref{DegreeThree}) below.

\pagebreak

\section{Discriminants of symmetric polynomials}

\subsection{Parametrisation of symmetric polynomials}

To derive a formula for discriminants of symmetric polynomials, let us recall that any symmetric polynomial in variables $x_1, \ldots, x_n$ can be expressed as a polynomial of elementary symmetric variables

\begin{align}
p_k = x_1^k + \ldots + x_n^k
\end{align}
\smallskip\\
For example, a homogeneous symmetric polynomial of degree $r = 2$ can be written as

\begin{align}
S(x_1, \ldots, x_n) = C_{2} p_2 + C_{11} p_1^2
\end{align}
\smallskip\\
where $C_{2}, C_{11}$ are two arbitrary parameters. For degree $r = 3$ such polynomial takes form

\begin{align}
S(x_1, \ldots, x_n) = C_{3} p_3 + C_{21} p_2 p_1 + C_{111} p_1^3
\end{align}
\smallskip\\
and contains three parameters $C_{3}, C_{21}, C_{111}$. Similarly, for arbitrary degree $r$ we have

\begin{align}
S(x_1, \ldots, x_n) = \sum\limits_{|Y| = r} C_Y \ \Big( p_{Y_1} p_{Y_2} \ldots \Big)
\label{GenSymm}
\end{align}
\smallskip\\
where the sum goes over all sequences $Y = (Y_1, Y_2, \ldots)$ with non-increasing elements $Y_1 \geq Y_2 \ldots \geq 0$ and fixed degree $|Y| = Y_1 + Y_2 + \ldots = r$. Such non-increasing sequences are known as partitions or Young diagrams. The number of free parameters is equal to the number of partitions of degree $r$, which we denote as $P(r)$. For example, $P(2) = 2$ corresponds to partitions $(2)$ and $(1,1)$. Similarly $P(3) = 3$ corresponds to partitions $(3), (2,1)$ and $(1,1,1)$. Several first values of $P(r)$ are

\begin{align}
\begin{array}{c|cccccccccccccc}
r & 1 & 2 & 3 & 4 & 5 & 6 & 7 & 8 & 9 & 10 \\
\\
P(r) & 1 & 2 & 3 & 5 & 7 & 11 & 15 & 22 & 30 & 42 \\
\end{array}
\end{align}
\smallskip\\
A symmetric homogeneous polynomial of degree $r$ has $P(r)$ independent coefficients $C_Y$. Its discriminant is a function of these coefficients, which we denote ${\cal D}_{n|r} \left( C \right)$. Note, that the number of free parameters of the polynomial does not depend on the number of variables $n$ (for fixed degree $r$ and for $n \geq r$). This important property of symmetric polynomials makes it possible to study the large-$n$ asymptotics of ${\cal D}_{n|r} \left( C \right)$, which is potentially interesting from the physical point of view.

Note also, that only the first $n$ variables $p_k$ are indepedent. Because of this, for $n < r$ the number of independent parameters of the polynomial is accidentally lower than $P(r)$. We do not consider these low-dimensional cases at all; from now on we assume that $n \geq r$. Our final expressions, however, are valid for $n < r$ as well (we state this without a proof and check for particular examples).

\subsection{Description of critical points}

Let us describe the critical points of (\ref{GenSymm}). As usual, to find the critical points of a polynomial, one needs to study the system of polynomial's derivatives. If $S$ is a homogeneous symmetric polynomial (\ref{GenSymm}), then its derivatives w.r.t. variables $x_i$ are, of course, not symmetric polynomials. Instead, they have a form

\begin{align}
\dfrac{\partial S}{\partial x_i} = \sum\limits_{k = 0}^{r-1} W_k(x_1, \ldots, x_n) \ x_i^{r-1-k}
\label{Derivatives}
\end{align}
\smallskip\\
or, what is the same,

\begin{align}
\left(
\begin{array}{cccc}
\partial S / \partial x_1\\
\\
\partial S / \partial x_2\\
\\
\ldots \\
\\
\partial S / \partial x_n\\
\end{array}\right) =
\left( \begin{array}{ccccccc}
1 & x_1 & x_1^2 & \ldots & x_1^{r-1} \\
\\
1 & x_2 & x_2^2 & \ldots & x_2^{r-1} \\
\\
\ldots\\
\\
1 & x_n & x_n^2 & \ldots & x_n^{r-1} \\
\\
\end{array} \right) \cdot
\left( \begin{array}{cccc}
W_{r-1}\\
\\
W_{r-2}\\
\\
\ldots \\
\\
W_{0}\\
\end{array} \right)
\label{Decompos}
\end{align}
\smallskip\\
where $W_k(x_1, \ldots, x_n)$ are are homogeneous symmetric polynomials of degree $k$:

\begin{align}
\begin{array}{lllll}
W_0 = r C_r\\
\\
W_1 = (r-1) C_{r-1,1} p_1\\
\\
W_2 = (r-2) \big( C_{r-2,2} p_2 + C_{r-2,1,1} p_1^2 \big)\\
\end{array}
\end{align}
\smallskip\\
and so on. It is convenient to consider two cases, $C_r = 0$ and $C_r \neq 0$. As a simple consequence of (\ref{Decompos}), in the first case the discriminant ${\cal D}_{n|r} \left( C \right)$ vanishes. Indeed, in this case the polynomial has a critical point

\begin{align}
(x_1, \ldots, x_n) = ( \omega^0, \omega^1, \ldots, \omega^{r-1}, 0, \ldots, 0), \ \ \ \ \ \omega = \exp \left( \dfrac{2 \pi i}{r} \right)
\label{NullVector}
\end{align}
\smallskip\\
For such a vector, all the variables $p_k$ with $0 < k < r$ vanish, therefore all $W_k(x_1, \ldots, x_n)$ with $0 < k < r$ vanish. The vanishing of $W_0(x_1, \ldots, x_n)$ follows from $C_r = 0$. As one can see, in this case all the derivatives (\ref{Decompos}) vanish together in a point (\ref{NullVector}), thus the discriminant of $S$ vanishes.

Let us consider the second case, when $C_r \neq 0$. We are going to show, that if vector $(x_1, \ldots, x_n)$ is a critical point, then there are at most $(r-1)$ distinct variables among $x_1, \ldots, x_n$. I.e, let us show that any critical point of the system (\ref{Decompos}) either has a form

\begin{equation}
\addtolength{\fboxsep}{5pt}
\boxed{
\begin{gathered}
(x_1, \ldots, x_n) = (\underbrace{y_1, \ldots, y_1}_{M_1}, \underbrace{y_2, \ldots, y_2}_{M_2}, \ldots, \underbrace{y_{r-1}, \ldots, y_{r-1}}_{M_{r-1}})
\end{gathered}
}\label{Distinct}
\end{equation}
\smallskip\\
where $M_1 + \ldots + M_{r-1} = n$ and $y_i \neq y_j$ for $i \neq j$, or is obtained from (\ref{Distinct}) by permutations of components. This is quite easy to show: indeed, in this case the system of equations under consideration

\begin{align}
\left( \begin{array}{ccccccc}
1 & x_1 & x_1^2 & \ldots & x_1^{r-1} \\
\\
1 & x_2 & x_2^2 & \ldots & x_2^{r-1} \\
\\
\ldots\\
\\
1 & x_n & x_n^2 & \ldots & x_n^{r-1} \\
\\
\end{array} \right) \cdot
\left( \begin{array}{cccc}
W_{r-1}\\
\\
W_{r-2}\\
\\
\ldots \\
\\
W_{0}\\
\end{array} \right) = \left(
\begin{array}{cccc}
0\\
\\
0\\
\\
\ldots \\
\\
0\\
\end{array}\right)
\label{DecomposSystem}
\end{align}
has non-zero solutions if and only if all $ r \times r$ minors of the rectangular Vandermonde matrix in the left hand side vanish. Different $r \times r$ minors of this $n \times r$ matrix are in one-to-one correspondence with different ways to choose $r$ variables out of $x_1, \ldots, x_n$. For each particular choice, the minor vanishes if and only if two equal variables are chosen. Consequently, the system (\ref{DecomposSystem}) has non-zero solutions if and only if every set of $r$ variables out of $x_1, \ldots, x_n$ contains a pair of equal variables. This implies, that there are at most $(r-1)$ distinct variables among $x_1, \ldots, x_n$. Vectors (\ref{Distinct}) will be called critical points of type $(M_1, \ldots, M_{r-1})$.

Note that, according to the permutation symmetry of the problem, existence of a critical point of type (\ref{Distinct}) implies existence of many other critical points, obtained by permutations of the components of (\ref{Distinct}). For example, let $r = 3$, $n = 4$ and $(M_1, M_2) = (1,3)$. If the vector $(y_1, y_2, y_2, y_2)$ is a critical point, then vectors $(y_2, y_1, y_2, y_2)$, $(y_2, y_2, y_1, y_2)$ and $(y_2, y_2, y_2, y_1)$ are critical points as well. In this case, there are four of them. The total number of critical points of type $(M_1, \ldots, M_{r-1})$ is obviously given by the multinomial coefficient $(M_1 + \ldots + M_{r-1})! / M_{1}! \ldots / M_{r-1}!$.

\subsection{Factorisation of the discriminant}

To summarize the above arguments, the polynomial $S$ is degenerate (has non-zero critical points) if and only if \emph{at least one of} vectors (\ref{Distinct}) is a solution of (\ref{DecomposSystem}). If this is the case, then the degree of degeneracy (the number of critical points) is $(M_1 + \ldots + M_{r-1})! / M_{1}! \ldots / M_{r-1}!$. Consequently, ${\cal D}_{n|r} \left( C \right)$ is a product of several factors, which are labeled by numbers $(M_1, \ldots, M_{r-1})$ and have appropriate multiplicities:

\begin{equation}
\addtolength{\fboxsep}{5pt}
\boxed{
\begin{gathered}
{\cal D}_{n|r} \left( C \right) = \alpha_n \ C_r^{\beta_{n|r}} \ \prod\limits_{M_1 + \ldots + M_{r-1} = n} \Big( d_{M}\left( C \right) \Big)^{ \dfrac{\#_M!}{(r-1)!} \dfrac{(M_1 + \ldots + M_{r-1})!}{M_1! \ldots M_{r-1}!}  }
\end{gathered}
}\label{DiscrimDecomp}
\end{equation}
\smallskip\\
where $\#_M$ is the number of zeroes among $M_1, \ldots, M_{r-1}$, \ $\alpha_n$ is a $C$-independent overall normalisation needed to satisfy eq. (\ref{Normalisation}), degree $\beta_{n|r}$ is fixed by the total degree of the discriminant

\begin{align}
\beta_{n|r} = \deg {\cal D}_{n|r} - \sum\limits_{M_1 + \ldots + M_{r-1} = n} \dfrac{\#_M! \deg d_{M}}{(r-1)!} \dfrac{(M_1 + \ldots + M_{r-1})!}{M_1! \ldots M_{r-1}!}
\label{BetaDegree}
\end{align}
and $d_{M}\left( C \right)$ is an irreducible homogeneous polynomial which vanishes if and only if (\ref{DecomposSystem}) has solutions of type (\ref{Distinct}). The multiplicity of each factor is multplied by $\#_M!/(r-1)!$ to avoid counting twice configurations which differ only by renaming of variables $y_i$. For example, consider vectors $(y_1, y_2, y_2, y_2)$ and $(y_2, y_1, y_1, y_1)$ which have type $(1,3)$ and $(3,1)$, respectively. Formally different, they differ only by renaming of variables $(y_1, y_2)$ and actually describe one and the same critical point. To handle this combinatorial problem we take a product over all decompositions $M_1 + \ldots + M_{r-1} = n$ without any ordering prescriptions, but divide by $(r-1)!$ -- the number of all permutations of $(r-1)$ letters -- and multiply by $\#_M!$ because zeroes are indistinguishable. \clearpage
\emph{}
\vspace{3ex}

\subsection{An explicit formula for $d_{M}\left( C \right)$ }

Decomposition (\ref{DiscrimDecomp}) is arguably the most important property of discriminants of symmetric polynomials. To finish our presentation, we give an explicit formula for $d_{M}\left( C \right)$. After substitution of

\begin{align}
(x_1, \ldots, x_n) = \Big(\underbrace{y_1, \ldots, y_1}_{M_1}, \underbrace{y_2, \ldots, y_2}_{M_2}, \ldots, \underbrace{y_{r-1}, \ldots, y_{r-1}}_{M_{r-1}}\Big)
\end{align}
\smallskip\\
the system of derivatives takes form

\begin{align}
\left(\dfrac{\partial S}{\partial x_1}, \ldots, \dfrac{\partial S}{\partial x_n}\right) = \Big(\underbrace{P^{(1)}_M, \ldots, P^{(1)}_M}_{M_1}, \underbrace{P^{(2)}_M, \ldots, P^{(2)}_M}_{M_2}, \ldots, \underbrace{P^{(r-1)}_M, \ldots, P^{(r-1)}_M}_{M_{r-1}}\Big)
\label{Y-Derivatives}
\end{align}
\smallskip\\
where $P^{(i)}_{M}$ are homogeneous of degree $r - 1$ in variables $y_i$. The number of variables $y_i$ is $(r-1)$ if all $M_i$ are positive, and can be less than $(r-1)$ if some $M_i = 0$. As follows from the permutation symmetry of $S(x_1, \ldots, x_n)$, the differences of derivatives are divisible by the corresponding differences of variables:

\begin{align} P^{(i)}_{M} - P^{(j)}_{M} = (y_i - y_j) P^{(ij)}_{M} \end{align}
\smallskip\\
where $P^{(ij)}_{M}$ are homogeneous polynomials of degree $r - 2$. Similarly, we find

\begin{align} y_i P^{(jk)}_{M} - y_j P^{(ik)}_{M} + y_k P^{(ij)}_{M} = (y_i - y_j) (y_i - y_k) (y_j - y_k) P^{(ijk)}_{M}
\end{align}
\smallskip\\
where $P^{(ijk)}_{M}$ are homogeneous polynomials of degree $r - 3$. This procedure of division may be continued further. It may be slightly more convenient to express the results in determinantal form:

\begin{align}
P^{(ij)}_{M} = \dfrac{ \det\limits_{2 \times 2}\left( \begin{array}{cc} 1 & P^{(i)}_M \\ \\ 1 & P^{(j)}_M \end{array} \right) }{ \det\limits_{2 \times 2} \left( \begin{array}{cc} 1 & y_i \\ \\ 1 & y_j \end{array} \right) } \ , \ \ \ \ \ \ \ P^{(ijk)}_{M} = \dfrac{ \det\limits_{3 \times 3} \left( \begin{array}{ccc} 1 & y_i & P^{(i)}_M \\ \\ 1 & y_j & P^{(j)}_M \\ \\ 1 & y_k & P^{(k)}_M \end{array} \right) }{ \det\limits_{3 \times 3} \left( \begin{array}{ccc} 1 & y_i & y_i^2 \\ \\ 1 & y_j & y_j^2 \\ \\ 1 & y_k & y_k^2 \end{array} \right) } \ , \ \ \  \ \ \ \ \ldots
\label{2and3levelquantities}
\end{align}
\smallskip\\
and so on. The results of division are polynomials in variables $y_i$, because $S(x_1, \ldots, x_n)$ is symmetric under permutations of $x_1, \ldots, x_n$. Generally for $1 \leq k \leq (r-1)$ we obtain homogeneous polynomials

\begin{align}
P^{(i_1 \ldots i_k)}_{M} = \dfrac{ \det\limits_{k \times k} \left( \begin{array}{ccccc} 1 & y_{i_1} & \ldots & y_{i_1}^{k-2} & P^{(i_1)}_M \\ \\ 1 & y_{i_2} & \ldots & y_{i_2}^{k-2} & P^{(i_2)}_M \\ \\ \ldots \\ \\ 1 & y_{i_k} & \ldots & y_{i_k}^{k-2} & P^{(i_k)}_M \\ \\\end{array} \right) }{ \det\limits_{k \times k} \left( \begin{array}{ccccc} 1 & y_{i_1} & \ldots & y_{i_1}^{k-2} & y_{i_1}^{k-1} \\ \\ 1 & y_{i_2} & \ldots & y_{i_2}^{k-2} & y_{i_2}^{k-1} \\ \\ \ldots \\ \\ 1 & y_{i_k} & \ldots & y_{i_k}^{k-2} & y_{i_k}^{k-1} \\ \\\end{array} \right) }
\end{align}
\smallskip\\
of degree $r - k$ in variables $y_i$. They are defined only for pairwise distinct upper indices $i_1, \ldots, i_k$ and are symmetric in them. For $k \geq r$ they are undefined: one can not choose $r$ distinct items out of $(r-1)$.

\paragraph{Proposition.} The irreducible factor $d_{M}\left( C \right)$ from (\ref{DiscrimDecomp}) is equal to the following resultant:

\begin{equation}
\addtolength{\fboxsep}{5pt}
\boxed{
\begin{gathered}
d_{M}\left( C \right) = R \left\{ \ \sum\limits_{i} P^{(i)}_M, \ \ \sum\limits_{i < j} P^{(i j)}_M, \ \ \sum\limits_{i < j < k} P^{(i j k)}_M, \ \ \ldots \ \right\}
\end{gathered}
}\label{MainFormula}
\end{equation}
\smallskip\\
where the system of equations in the right hand side contains exactly as many equations as there are $y$-variables (this number is $(r-1)$ in general position and can be less than $(r-1)$ if some $M_i$ are vanishing).

\paragraph{Sketch of a proof.} It is easy to see, that any solution of $P^{(i)}_{M} = 0$ is a solution of $P^{(i_1, \ldots, i_k)}_{M} = 0$. Therefore, the resultant in the right hand side of (\ref{MainFormula}) must be divisible on $d_{M}\left( C \right)$:

\begin{align}
R \left\{ \ \sum\limits_{i} P^{(i)}_M, \ \ \sum\limits_{i < j} P^{(i j)}_M, \ \ \sum\limits_{i < j < k} P^{(i j k)}_M, \ \ \ldots \ \right\} = d_{M}\left( C \right) {\widetilde d}_{M}\left( C \right)
\label{Factors2}
\end{align}
\smallskip\\
The unknown factor ${\widetilde d}_{M}\left( C \right)$ can, in principle, be some polynomial in coefficients $C_Y$.
To complete the proof, one just needs to prove that the resultant on the left hand side of (\ref{Factors2}) is an irreducible polynomial in $C_Y$. Then, the other factor must have degree zero in $C_Y$ and is nothing but a constant of proportionality (which can be put to unity, since (\ref{MainFormula}) already contains a normalisation constant $\alpha_n$).

Unfortunately, we are yet unable to prove irreducibility of this resultant for arbitrary $r$. However, in particular examples -- for $r = 2,3,4$ -- this resultant can be computed and turns out to be irreducible. This justifies the use of (\ref{MainFormula}) in these examples. The proof for arbitrary $r$ remains to be done.
\smallskip\\

Together, eqs. (\ref{DiscrimDecomp}) and (\ref{MainFormula}) constitute our main result - an explicit and division free formula for discriminants of symmetric homogeneous polynomials of degree $r$ in $n \geq r$ variables, expressing them through resultants in no more than $(r-1)$ variables. As we show below, this formula allows to calculate such discriminants in dimensions $n >> r$, where other algorithms become non-practical.

\section{Symmetric polynomials of degree 2}

We begin our presentation of examples with the case of symmetric polynomials of degree $2$, i.e, the case of quadratic forms. This example is quite simple, since it is completely treatable by methods of linear algebra. A homogeneous symmetric polynomial of degree two should have a form

\begin{align}
S(x_1, \ldots, x_n) = C_{2} p_2 + C_{11} p_1^2
\end{align}
\smallskip\\
where $C_{2}, C_{11}$ are two arbitrary parameters. To calculate the discriminant, we take derivatives and obtain

\begin{align}
\dfrac{\partial S}{\partial x_i} = 2 C_{2} x_i + 2 C_{11} p_1 = 0
\label{Sys2}
\end{align}
\smallskip\\
The simplest option is just to write a matrix and calculate its determinant by usual rules:

\begin{align}
{\cal D}_{n|2}\big(C_{2}, C_{11}\big) = \det\limits_{n \times n} \left(
\begin{array}{cccc}
C_{2} + C_{11} & C_{11} & \ldots & C_{11} \\
\\
C_{11} & C_{2} + C_{11} & \ldots & C_{11} \\
\\
\ldots \\
C_{11} & C_{11} & \ldots & C_{2} + C_{11} \\
\end{array} \right) = C_{2}^{n-1} \big( C_2 + n C_{11} \big)
\end{align}
\smallskip\\
As one can see, the discriminant is highly factorized. As we know from the previous sections, this is a general property  of symmetric polynomials. It is this property that makes discriminants of symmetric polynomials such an interesting and simple object of investigation. Now let us reproduce the same answer from eq. (\ref{DiscrimDecomp}), which in this case takes form

\begin{align}
{\cal D}_{n|2}\big(C_{2}, C_{11}\big) = \alpha_n C_2^{\beta_{n|2}} d_{n}\big(C_{2}, C_{11}\big)
\end{align}
\smallskip\\
i.e, in this case we have a single decomposition $M = (n)$ of the number $n$ into $(r-1) = 1$ parts. Accordingly, there is a single polynomial $d_{n}\left( C \right)$ in this case. Direct calculation with (\ref{Y-Derivatives}) gives

\begin{align}
P^{(1)}_n(y) = \big( 2 n C_{2} + 2 n^2 C_{11} \big) y
\end{align}
\smallskip\\
Note, that in this very simple example the resultant is taken in a single variable $y$:

\begin{align}
d_{n}\left( C \right) = R \left\{ \sum\limits_{i} P^{(i)}_n \right\} = R \left\{ 2 n C_{2} y + 2 n^2 C_{11} y \right\} = 2 n C_{2} + 2 n^2 C_{11}
\end{align}
\smallskip\\
This is because, as follows from (\ref{Distinct}), permutation symmetry forces any critical point to have a form $(y, \ldots, y)$. It is clear now that $\deg d_{n}\left( C \right) = 1$, so that $\beta_{n|2} = \deg {\cal D}_{n|2} - \deg d_{n} = n - 1$ and

\begin{align}
{\cal D}_{n|2}\big(C_{2}, C_{11}\big) = 2 n \alpha_n C_2^{n-1} \big( C_{2} + n C_{11} \big)
\end{align}
\smallskip\\
Fixing as in eq. (\ref{Normalisation}) the normalisation, we find

\begin{align}
{\cal D}_{n|2}\big(1, 0\big) = 2n \alpha_n = 1, \ \ \ \ \ \alpha_n = \dfrac{1}{2n}
\end{align}
\smallskip\\
and, finally,

\begin{equation}
\addtolength{\fboxsep}{5pt}
\boxed{
\begin{gathered}
{\cal D}_{n|2}\big(C_{2}, C_{11}\big) = C_2^{n-1} \big( C_{2} + n C_{11} \big)
\end{gathered}
}\label{DegreeTwo}
\end{equation}
\smallskip\\
This is an explicit formula for the discriminant of a symmetric quadratic form in arbitrarily high dimension $n$.
This basic example can be considered as a simple illustration of what happens in the general situation.

\section{Symmetric polynomials of degree 3}

A homogeneous symmetric polynomial of degree $3$ in $n$ variables has a form

\begin{align}
S(x_1, \ldots, x_n) = C_{3} p_3 + C_{21} p_2 p_1 + C_{111} p_1^3
\end{align}
\smallskip\\
and contains three parameters $C_{3}, C_{21}, C_{111}$. To calculate the discriminant, we take derivatives and obtain

\begin{align}
\dfrac{\partial S}{\partial x_i} = 3 C_{3} x_i^2 + 2 C_{21} p_1 x_i + C_{21} p_2 + 3 C_{111} p_1^2
\end{align}
\smallskip\\
The main formula (\ref{DiscrimDecomp}) in this case turns into a product over all decompositions of $n$ into two parts:

\begin{align}
{\cal D}_{n|3}\big(C_{3}, C_{21}, C_{111}\big) = \alpha_n C_3^{\beta_{n|3}} \prod\limits_{M_1 + M_2 = n} \Big( d_{M_1 M_2}\big(C_{3}, C_{21}, C_{111}\big) \Big)^{\dfrac{n!}{2 M_1!M_2!} }
\label{3Main}
\end{align}
\smallskip\\
(where $\#_M$ equals either 0 or 1, therefore $\#_M! = 1$). We now need to calculate $d_{M_1 M_2}\big(C_{3}, C_{21}, C_{111}\big)$. First, let us consider the case $M_1, M_2 \neq 0$. Following section 2.4, we express the derivatives through the variables $y_i$ and the parameters $M_i$:

{\fontsize{9pt}{0pt}
\begin{align}
P^{(1)}_{M_1, M_2}(y_1, y_2) = 3( C_{3} + C_{21} M_1 + C_{111} M_1^2 ) y_1^2 + ( 2 C_{21} M_2 + 6 C_{111} M_1 M_2 ) y_1 y_2 + (C_{21} M_2 + 3 C_{111} M_2^2 ) y_2^2\\
\nonumber \\
P^{(2)}_{M_1, M_2}(y_1, y_2) = (C_{21} M_1 + 3 C_{111} M_1^2 ) y_1^2 + ( 2 C_{21} M_1 + 6 C_{111} M_1 M_2 ) y_1 y_2 + 3 (C_{3} + C_{21} M_2 + C_{111} M_2^2 ) y_2^2
\end{align}}
\smallskip\\
Their sum is equal to

{\fontsize{9pt}{0pt}\begin{align}
P^{(1)}_{M_1, M_2}(y_1, y_2) + P^{(2)}_{M_1, M_2}(y_1, y_2) \ = \ & (4 C_{2 1} M_{1}+6 C_{1 1 1} M_{1}^2+3 C_{3}) y_{1}^2 + \emph{} \\ \nonumber & \\ \nonumber & \emph{} + (2 C_{2 1} M_{1}+12 C_{1 1 1} M_{1} M_{2}+2 C_{2 1} M_{2}) y_{1} y_{2} + \emph{} \\ \nonumber & \\ \nonumber & \emph{} + (3 C_{3}+4 C_{2 1} M_{2}+6 C_{1 1 1} M_{2}^2) y_{2}^2
\end{align}}
\smallskip\\
Using eq. (\ref{2and3levelquantities}), we then calculate the second-level quantity

{\fontsize{9pt}{0pt}\begin{align}
P^{(12)}_{M_1, M_2}(y_1, y_2) = \dfrac{P^{(1)}_{M_1,M_2} - P^{(2)}_{M_1,M_2}}{y_1 - y_2} = (3 C_{3}+2 C_{2 1} M_{1}) y_{1}+(2 C_{2 1} M_{2}+3 C_{3}) y_{2}
\end{align}}
\smallskip\\
According to (\ref{MainFormula}),  the irreducible factors are equal to

{\fontsize{9pt}{0pt}\begin{align}
d_{M_1 M_2}\big(C_{3}, C_{21}, C_{111}\big) = R\Big\{ P^{(1)}_{M_1, M_2} + P^{(2)}_{M_1, M_2}, P^{(12)}_{M_1, M_2} \Big\}
\end{align}}
\smallskip\\
The resultant here is taken in two variables $y_1, y_2$ and, therefore, is just a Sylvester resultant \cite{GKZ, NOLINAL}:

{\fontsize{9pt}{0pt}\begin{align*}
 d_{M_1 M_2}\big(C_{3}, C_{21}, C_{111}\big) =  R\Big\{ P^{(1)}_{M_1, M_2}(y_1, y_2) + P^{(2)}_{M_1, M_2}(y_1, y_2), P^{(12)}_{M_1, M_2}(y_1, y_2) \Big\} =
\end{align*}}

{\fontsize{9pt}{0pt}\begin{align*}
= \det\limits_{3 \times 3} \left(
\begin{array}{cccc}
3 C_{3}+4 C_{2 1} M_{2}+6 C_{1 1 1} M_{2}^2 & 2 C_{2 1} M_{1}+12 C_{1 1 1} M_{1} M_{2}+2 C_{2 1} M_{2} & 4 C_{2 1} M_{1}+6 C_{1 1 1} M_{1}^2+3 C_{3} \\
\\
2 C_{2 1} M_{2}+3 C_{3} & 3 C_{3}+2 C_{2 1} M_{1} & 0 \\
\\
0 & 2 C_{2 1} M_{2}+3 C_{3} & 3 C_{3}+2 C_{2 1} M_{1} \\
\\
\end{array}
\right) =
\end{align*}}

{\fontsize{9pt}{0pt}\begin{align}
\nonumber \ = \ & \big(-216 C_{1 1 1} C_{3}^2+72 C_{2 1}^2 C_{3}+8 C_{2 1}^3 (M_1 + M_2)\big) M_1 M_2 + \emph{} \\ \nonumber & \\ & \emph{} + 54 C_{3}^3+54 C_{2 1} C_{3}^2 (M_1+M_2)+54 C_{1 1 1} C_{3}^2 (M_1+M_2)^2
\end{align}}
\smallskip\\
Since $M_1 + M_2 = n$, for $M_1, M_2 \neq 0$ we finally obtain

{\fontsize{9pt}{0pt}\begin{align}
d_{M_1 M_2}\big(C_{3}, C_{21}, C_{111}\big) = \big(8 C_{2 1}^3 n + 72 C_{2 1}^2 C_{3} - 216 C_{1 1 1} C_{3}^2 \big) M_1 M_2 + 54 \big( C_{3}^3 + C_{2 1} C_{3}^2 n + C_{1 1 1} C_{3}^2 n^2 \big)
\label{X1}
\end{align}}
\smallskip\\
If $M_1 = 0$ or $M_2 = 0$, things get even simpler. In this case the polynomials $P$ depend on a single variable

{\fontsize{9pt}{0pt}\begin{align}
P^{(1)}_{n,0}(y_1) = 3 ( C_2 + C_{21} n + C_{111} n^2 ) y_1^2, \ \ \ \ \ P^{(1)}_{0,n}(y_2) = 3 ( C_2 + C_{21} n + C_{111} n^2 ) y_2^2
\end{align}}
\smallskip\\
and, accordingly, resultants are taken in a single variable:

{\fontsize{9pt}{0pt}\begin{align}
d_{n 0}\big(C_{3}, C_{21}, C_{111}\big) = d_{0 n}\big(C_{3}, C_{21}, C_{111}\big) = 3 ( C_2 + C_{21} n + C_{111} n^2 )
\label{X2}
\end{align}}
\smallskip\\
Note, that this expression differs by a factor of $18 C_3^2$ from the formal limit $M_1 = 0$ or $M_2 = 0$ of (\ref{X1}). Therefore, it is essential to consider the cases with vanishing $M_i$ separately. Substituting (\ref{X1}) and (\ref{X2}) into (\ref{3Main}) we find after some simplifications

{\fontsize{9pt}{0pt}\begin{align*}
{\cal D}_{n|3}\big(C_{3}, C_{21}, C_{111}\big) =
\end{align*}}
{\fontsize{9pt}{0pt}
\begin{align*}
= \dfrac{ 54^{2^{n-1}} \alpha_n}{18} \ C_3^{\beta_{n|3} - 2}  \prod\limits_{m = 0}^{n} \left[ 4 m(n-m) \left( \dfrac{C_{2 1}^3}{27} n + \dfrac{C_{2 1}^2 C_{3}}{3} - C_{1 1 1} C_{3}^2 \right) + C_{3}^2 \big( C_3 + C_{2 1} n + C_{1 1 1} n^2 \big) \right]^{\dfrac{n!}{2 m!(n-m)!}}
\end{align*}}
\smallskip\\
Calculating with the help of eq. (\ref{BetaDegree}) the degree

{\fontsize{9pt}{0pt}\begin{align}
\beta_{n|3} = n 2^{n-1} - \sum\limits_{m = 1}^{n-1} \dfrac{3 n!}{2 m!(n-m)!} - 1 = (n - 3) 2^{n-1} + 2
\end{align}}
\smallskip\\
and fixing as in eq. (\ref{Normalisation}) the normalisation

{\fontsize{9pt}{0pt}\begin{align}
{\cal D}_{n|3}\big(1, 0, 0\big) = \dfrac{ 54^{2^{n-1}} \alpha_n}{18} = 1, \ \ \ \ \ \alpha_n = 18 \cdot 54^{-2^{n-1}}
\end{align}}
\smallskip\\
we obtain, finally,

{\fontsize{9pt}{0pt}
\begin{equation}
\addtolength{\fboxsep}{5pt}
\boxed{
\begin{gathered}
{\cal D}_{n|3}\big(C_{3}, C_{21}, C_{111}\big) = C_3^{(n - 3) 2^{n-1}} \times \emph{} \\
\\
\emph{} \times \prod\limits_{m = 0}^{n} \left[ 4 m(n-m) \left( \dfrac{1}{27} C_{2 1}^3 n + \dfrac{1}{3} C_{2 1}^2 C_{3} - C_{1 1 1} C_{3}^2 \right) + C_{3}^2 \big( C_3 + C_{2 1} n + C_{1 1 1} n^2 \big) \right]^{\dfrac{n!}{2 m!(n-m)!}}
\end{gathered}
}\label{DegreeThree}
\end{equation}}
\smallskip\\
This is an explicit formula for the discriminant of a symmetric cubic in arbitrarily high dimension $n \geq 3$. To demonstrate the computational power of this formula, let us put here, say, $n = 20$. As far as we know, most of the known algorithms become inpractical already for $n \approx 10$. A fast computer calculation gives

\begin{align}
\nonumber {\cal D}_{20|3}\big(C_{3}, C_{21}, C_{111}\big) \ = \ & C_{3}^{8912896} \ (400 C_{1 1 1} C_{3}^{2}+20 C_{2 1} C_{3}^2+C_{3}^3) \times \emph{} \\ \nonumber & \\ \nonumber & \emph{}
(324 C_{1 1 1} C_{3}^2+\dfrac{1520}{27} C_{2 1}^3 + \dfrac{76}{3} C_{2 1}^2 C_{3}+20 C_{2 1} C_{3}^2+C_{3}^3)^{20}
\times \emph{} \\ \nonumber & \\ \nonumber & \emph{}
(256 C_{1 1 1} C_{3}^2+\dfrac{320}{3} C_{2 1}^3+48 C_{2 1}^2 C_{3}+20 C_{2 1} C_{3}^2+C_{3}^3)^{190}
\times \emph{} \\ \nonumber & \\ \nonumber & \emph{}
(196 C_{1 1 1} C_{3}^2+\dfrac{1360}{9} C_{2 1}^3+68 C_{2 1}^2 C_{3}+20 C_{2 1} C_{3}^2+C_{3}^3)^{1140}
\times \emph{} \\ \nonumber & \\ \nonumber & \emph{}
(144 C_{1 1 1} C_{3}^2+\dfrac{5120}{27} C_{2 1}^3+\dfrac{256}{3} C_{2 1}^2 C_{3}+20 C_{2 1} C_{3}^2+C_{3}^3)^{4845}
\times \emph{} \\ \nonumber & \\ \nonumber & \emph{}
(100 C_{1 1 1} C_{3}^2+\dfrac{2000}{9} C_{2 1}^3+100 C_{2 1}^2 C_{3}+20 C_{2 1} C_{3}^2+C_{3}^3)^{15504}
\times \emph{} \\ \nonumber & \\ \nonumber & \emph{}
(64 C_{1 1 1} C_{3}^2+\dfrac{2240}{9} C_{2 1}^3+112 C_{2 1}^2 C_{3}+20 C_{2 1} C_{3}^2+C_{3}^3)^{38760}
\times \emph{} \\ \nonumber & \\ \nonumber & \emph{}
(36 C_{1 1 1} C_{3}^2+\dfrac{7280}{27} C_{2 1}^3+\dfrac{364}{3} C_{2 1}^2 C_{3}+20 C_{2 1} C_{3}^2+C_{3}^3)^{77520}
\times \emph{} \\ \nonumber & \\ \nonumber & \emph{}
(16 C_{1 1 1} C_{3}^2+\dfrac{2560}{9} C_{2 1}^3+128 C_{2 1}^2 C_{3}+20 C_{2 1} C_{3}^2+C_{3}^3)^{125970}
\times \emph{} \\ \nonumber & \\ \nonumber & \emph{}
(4 C_{1 1 1} C_{3}^2+\dfrac{880}{3} C_{2 1}^3+132 C_{2 1}^2 C_{3}+20 C_{2 1} C_{3}^2+C_{3}^3)^{167960}
\times \emph{} \\ \nonumber & \\ & \emph{}
(\dfrac{8000}{27} C_{2 1}^3+\dfrac{400}{3} C_{2 1}^2 C_{3}+20 C_{2 1} C_{3}^2+C_{3}^3)^{92378}
\end{align}
\smallskip\\
This expression is somewhat charming -- it is a closed formula for a certain discriminant in 20 variables. It is permutation symmetry of $S$ that allows to write such closed formulas. Note the degree of $C_3$, which is equal to $8912896 = 17 \cdot 2^{19}$.  Note also, that even a slight deformation of the polynomial $S$, which destroys the permutation symmetry ( say, adding a term $+ \epsilon x_1 x_{2} x_{3}$ ) destroys also this strong factorisation of the discriminant and results in billions of billions of terms, which are hard to imagine.

For small values of $n$, an independent check of eq. (\ref{DegreeThree}) can be made. For example, for arbitrary (not necessarily symmetric) homogeneous cubic polynomial $S(x_1,x_2,x_3)$ in three variables

\begin{align}
\nonumber S(x_1,x_2,x_3) \ = \ & S_{1 1 1} x_{1}^3+3 S_{1 1 2} x_{1}^2 x_{2} + 3 S_{1 1 3} x_{1}^2 x_{3} + 3 S_{1 2 2} x_{1} x_{2}^2 + 6 S_{1 2 3} x_{1} x_{2} x_{3} + \emph{} \\ \nonumber & \\ & \emph{} + 3 S_{1 3 3} x_{1} x_{3}^2 + S_{2 2 2} x_{2}^3 + 3 S_{2 2 3} x_{2}^2 x_{3} + 3 S_{2 3 3} x_{2} x_{3}^2 + S_{3 3 3} x_{3}^3
\end{align}
\smallskip\\
there exists a well-known explicit formula for the discriminant: \clearpage

\emph{}
\vspace{3ex}

\begin{align}
D_{3|3}(S) = \det\limits_{6 \times 6} \left( \begin{array}{ccccccccc}
S_{111} & S_{112} & S_{113} & S_{122} & S_{123} & S_{133} \\ \\
S_{112} & S_{122} & S_{123} & S_{222} & S_{223} & S_{233} \\ \\
S_{113} & S_{123} & S_{133} & S_{223} & S_{233} & S_{333} \\ \\
H_{111} & H_{112} & H_{113} & H_{122} & H_{123} & H_{133} \\ \\
H_{112} & H_{122} & H_{123} & H_{222} & H_{223} & H_{233} \\ \\
H_{113} & H_{123} & H_{133} & H_{223} & H_{233} & H_{333}
\end{array} \right)
\label{Sylv3}
\end{align}
\smallskip\\
where $H_{ijk}$ are expansion coefficients of the Hessian determinant

\begin{align}
\nonumber H(x_1,x_2,x_3) = \det\limits_{3 \times 3} \left( \dfrac{\partial^2 S}{\partial x_i \partial x_j} \right) \ = \ & H_{1 1 1} x_{1}^3+3 H_{1 1 2} x_{1}^2 x_{2} + 3 H_{1 1 3} x_{1}^2 x_{3} + 3 H_{1 2 2} x_{1} x_{2}^2 + 6 H_{1 2 3} x_{1} x_{2} x_{3} + \emph{} \\ \nonumber & \\ & \emph{} + 3 H_{1 3 3} x_{1} x_{3}^2 + H_{2 2 2} x_{2}^3 + 3 H_{2 2 3} x_{2}^2 x_{3} + 3 H_{2 3 3} x_{2} x_{3}^2 + H_{3 3 3} x_{3}^3
\end{align}
\smallskip\\
Formula (\ref{Sylv3}) was first found by Sylvester, as a generalisation of his well-known formula for $n = 2$ resultants and discriminants. In the case of symmetric polynomial

\begin{align}
S(x_1,x_2,x_3) = C_3 (x_1^3 + x_2^3 + x_3^3) + C_{21} (x_1^2 + x_2^2 + x_3^2) (x_1 + x_2 + x_3) + C_{111} (x_1 + x_2 + x_3)^3
\end{align}
\smallskip\\
our formula gives

\begin{align}
{\cal D}_{3|3}\big(C_{3}, C_{21}, C_{111}\big) = (9 C_{1 1 1} C_{3}^2+3 C_{2 1} C_{3}^2+C_{3}^3) (C_{1 1 1} C_{3}^2 + \dfrac{8}{9} C_{2 1}^3 + \dfrac{8}{3} C_{2 1}^2 C_{3}+3 C_{2 1} C_{3}^2+C_{3}^3)^3
\end{align}
\smallskip\\
and one can easily check that calculation with the $6 \times 6$ matrix gives the same result. This can be regarded as an independent check of (\ref{DegreeThree}). A similar check can be made for $n = 4$ with a $20 \times 20$ matrix, which is analogous to (\ref{Sylv3}) and which we do not include here. For $n > 4$ a check of this type becomes impossible, since no analogues of (\ref{Sylv3}) are known for $n > 4$.

The formula (\ref{DiscrimDecomp}) was derived in the assumption of $n \geq r$. To check its validity in the region $n < r$, let us consider now the case $n = 2$. For arbitrary (not necessarily symmetric) homogeneous cubic polynomial $S(x_1,x_2)$ in two variables

\begin{align}
S(x_1,x_2) = S_{111} x_1^3 + 3 S_{112} x_1^2 x_2 + 3 S_{122} x_1 x_2^2 + S_{222} x_2^3
\end{align}
\smallskip\\
the discriminant is easily calculated with ordinary Sylvester matrix and equals

\begin{align}
D_{2|3}(S) & \ = \ \det\limits_{4 \times 4} \left( \begin{array}{ccccc}
S_{111} & 2 S_{112} & S_{122} & 0 \\ \\
0 & S_{111} & 2 S_{112} & S_{122} \\ \\
S_{112} & 2 S_{122} & S_{222} & 0 \\ \\
0 & S_{112} & 2 S_{122} & S_{222} \end{array} \right) = \emph{} \\ \nonumber & \\ & \emph{} = S_{1 1 1}^2 S_{2 2 2}^2-6 S_{1 1 1} S_{1 1 2} S_{1 2 2} S_{2 2 2}+4 S_{1 1 1} S_{1 2 2}^3+4 S_{1 1 2}^3 S_{2 2 2}-3 S_{1 1 2}^2 S_{1 2 2}^2
\label{Discrim23}
\end{align}
\smallskip\\
In the case of symmetric polynomial

\begin{align}
S(x_1,x_2) = C_3 (x_1^3 + x_2^3) + C_{21} (x_1^2 + x_2^2) (x_1 + x_2) + C_{111} (x_1 + x_2)^3
\end{align}
we find

\begin{align}
S_{111} = S_{222} = C_3 + C_{21} + C_{111}, \ \ \ \ \ S_{112} = S_{122} = C_{111} + \dfrac{1}{3} C_{21}
\end{align}
\smallskip\\
Substituting these coefficients into (\ref{Discrim23}) we find

\begin{align}
D_{2|3}\Big( S \Big) = \dfrac{1}{27} (4 C_{1 1 1}+2 C_{2 1}+C_{3}) (3 C_{3}+2 C_{2 1})^3
\end{align}
\smallskip\\
At the same time, our formula (\ref{DegreeThree}) gives

\begin{align}
{\cal D}_{2|3}\big(C_{3}, C_{21}, C_{111}\big) \ = \ \dfrac{1}{27} (4 C_{1 1 1}+2 C_{2 1}+C_{3}) (3 C_{3}+2 C_{2 1})^3
\end{align}
\smallskip\\
Therefore, (\ref{DegreeThree}) stays valid even for $n < r$. To conclude this section, we note that eq. (\ref{DegreeThree}) can be rewritten, after some algebraic transformations, in even more concise form:

{\fontsize{9pt}{0pt}
\begin{equation}
\addtolength{\fboxsep}{5pt}
\boxed{
\begin{gathered}
{\cal D}_{n|3}\big(C_{3}, C_{21}, C_{111}\big) = \big( B_3 \big)^{(n-3)2^{n-1}} \prod\limits_{k = 0}^{n - 1} \left( \left( \dfrac{n-2k}{9n} \right)^2 B_1 B_3^2 + \dfrac{4k(n-k)}{27 n^2} B_2^3 \right)^{\dfrac{(n-1)!}{k!(n-1-k)!}}
\end{gathered}
}\label{DegreeThree2}
\end{equation}}
\smallskip\\
where

\begin{align}
\left\{
\begin{array}{lll}
B_{1} = n^{2}C_{111} + n C_{21} + C_{3},\\
\\
B_{2} = n C_{21} + 3 C_{3},\\
\\
B_{3} = C_{3},
\end{array}
\right.
\end{align}
\smallskip\\
are just another parameters in the space of cubic symmetric polynomials. The present paper is devoted to generalisation of this beautiful formula to arbitrary degrees $r$.

\section{Symmetric polynomials of degree 4}

A homogeneous symmetric polynomial of degree $4$ in $n$ variables has a form

\begin{align}
S(x_1, \ldots, x_n) = C_{4} p_4 + C_{31} p_3 p_1 + C_{22} p_2^2 + C_{211} p_2 p_1^2 + C_{1111} p_1^4
\end{align}
\smallskip\\
and contains five parameters $C_{4}, C_{31}, C_{22}, C_{211}, C_{1111}$. To calculate the discriminant, we take derivatives

\begin{align}
\dfrac{\partial S}{\partial x_i} = 4 C_{4} x_i^3 + 3 C_{31} p_1 x_i^2 + (4 C_{22} p_2 + 2 C_{211} p_1^2) x_i + ( C_{31} p_3 + 2 C_{211} p_2 p_1 + 4 C_{1111} p_1^3 )
\end{align}
\smallskip\\
The main formula (\ref{DiscrimDecomp}) in this case turns into a product over all decompositions of $n$ into three parts:

\begin{align}
{\cal D}_{n|4}\big(C_{4}, C_{31}, C_{22}, C_{211}, C_{1111}\big) = \alpha_n C_4^{\beta_{n|4}} \prod\limits_{M_1 + M_2 + M_3 = n} \Big( d_{M_1 M_2 M_3}\big( C \big) \Big)^{\dfrac{\#_M! n!}{6 M_1!M_2!M_3!} }
\label{4Main}
\end{align}
\smallskip\\
We now need to calculate $d_{M_1 M_2 M_3}\big( C \big)$. First, let us consider the case $M_1, M_2, M_3 \neq 0$. Following section 2.4, we express the derivatives through the variables $y_i$ and the parameters $M_i$:

\begin{align}
\nonumber P^{(1)}_{M_1, M_2, M_3}(y_1, y_2, y_3) \ = \ &
(4 C_{2 2} M_{1}+4 C_{1 1 1 1} M_{1}^3+4 C_{3 1} M_{1}+4 C_{2 1 1} M_{1}^2+4 C_{4}) y_{1}^3 + \emph{} \\ \nonumber & \\ \nonumber & \emph{} +
(12 C_{1 1 1 1} M_{1}^2 M_{2}+3 C_{3 1} M_{2}+6 C_{2 1 1} M_{1} M_{2}) y_{1}^2 y_{2}+ \emph{} \\ \nonumber & \\ \nonumber & \emph{} +
(12 C_{1 1 1 1} M_{1}^2 M_{3}+3 C_{3 1} M_{3}+6 C_{2 1 1} M_{1} M_{3}) y_{1}^2 y_{3}+ \emph{} \\ \nonumber & \\ \nonumber & \emph{} +
(4 C_{2 2} M_{2}+2 C_{2 1 1} M_{1} M_{2}+2 C_{2 1 1} M_{2}^2+12 C_{1 1 1 1} M_{1} M_{2}^2) y_{1} y_{2}^2+ \emph{} \\ \nonumber & \\ \nonumber & \emph{} +
(4 C_{2 1 1} M_{2} M_{3}+24 C_{1 1 1 1} M_{1} M_{2} M_{3}) y_{1} y_{2} y_{3}+ \emph{} \\ \nonumber & \\ \nonumber & \emph{} +
(4 C_{2 2} M_{3}+2 C_{2 1 1} M_{1} M_{3}+12 C_{1 1 1 1} M_{1} M_{3}^2+2 C_{2 1 1} M_{3}^2) y_{1} y_{3}^2+ \emph{} \\ \nonumber & \\ \nonumber & \emph{} +
(4 C_{1 1 1 1} M_{2}^3+C_{3 1} M_{2}+2 C_{2 1 1} M_{2}^2) y_{2}^3+ \emph{} \\ \nonumber & \\ \nonumber & \emph{} +
(12 C_{1 1 1 1} M_{2}^2 M_{3}+2 C_{2 1 1} M_{2} M_{3}) y_{2}^2 y_{3}+ \emph{} \\ \nonumber & \\ \nonumber & \emph{} +
(2 C_{2 1 1} M_{2} M_{3}+12 C_{1 1 1 1} M_{2} M_{3}^2) y_{2} y_{3}^2+ \emph{} \\ \nonumber & \\ & \emph{} +
(C_{3 1} M_{3}+2 C_{2 1 1} M_{3}^2+4 C_{1 1 1 1} M_{3}^3) y_{3}^3
\end{align}

\begin{align}
\nonumber P^{(2)}_{M_1, M_2, M_3}(y_1, y_2, y_3) \ = \ &
(2 C_{2 1 1} M_{1}^2+4 C_{1 1 1 1} M_{1}^3+C_{3 1} M_{1}) y_{1}^3+ \emph{} \\ \nonumber & \\ \nonumber & \emph{} +(12 C_{1 1 1 1} M_{1}^2 M_{2}+2 C_{2 1 1} M_{1} M_{2}+4 C_{2 2} M_{1}+2 C_{2 1 1} M_{1}^2) y_{1}^2 y_{2}+ \emph{} \\ \nonumber & \\ \nonumber & \emph{} +(12 C_{1 1 1 1} M_{1}^2 M_{3}+2 C_{2 1 1} M_{1} M_{3}) y_{1}^2 y_{3}+ \emph{} \\ \nonumber & \\ \nonumber & \emph{} +(6 C_{2 1 1} M_{1} M_{2}+3 C_{3 1} M_{1}+12 C_{1 1 1 1} M_{1} M_{2}^2) y_{1} y_{2}^2+ \emph{} \\ \nonumber & \\ \nonumber & \emph{} +(4 C_{2 1 1} M_{1} M_{3}+24 C_{1 1 1 1} M_{1} M_{2} M_{3}) y_{1} y_{2} y_{3}+ \emph{} \\ \nonumber & \\ \nonumber & \emph{} +(12 C_{1 1 1 1} M_{1} M_{3}^2+2 C_{2 1 1} M_{1} M_{3}) y_{1} y_{3}^2+ \emph{} \\ \nonumber & \\ \nonumber & \emph{} +(4 C_{2 2} M_{2}+4 C_{1 1 1 1} M_{2}^3+4 C_{2 1 1} M_{2}^2+4 C_{3 1} M_{2}+4 C_{4}) y_{2}^3+ \emph{} \\ \nonumber & \\ \nonumber & \emph{} +(3 C_{3 1} M_{3}+6 C_{2 1 1} M_{2} M_{3}+12 C_{1 1 1 1} M_{2}^2 M_{3}) y_{2}^2 y_{3}+ \emph{} \\ \nonumber & \\ \nonumber & \emph{} +(2 C_{2 1 1} M_{2} M_{3}+2 C_{2 1 1} M_{3}^2+12 C_{1 1 1 1} M_{2} M_{3}^2+4 C_{2 2} M_{3}) y_{2} y_{3}^2+ \emph{} \\ \nonumber & \\ & \emph{} +(C_{3 1} M_{3}+2 C_{2 1 1} M_{3}^2+4 C_{1 1 1 1} M_{3}^3) y_{3}^3 \\ \nonumber
\\ \nonumber
P^{(3)}_{M_1, M_2, M_3}(y_1, y_2, y_3) \ = \ &
(2 C_{2 1 1} M_{1}^2+4 C_{1 1 1 1} M_{1}^3+C_{3 1} M_{1}) y_{1}^3+ \emph{} \\ \nonumber & \\ \nonumber & \emph{} +(12 C_{1 1 1 1} M_{1}^2 M_{2}+2 C_{2 1 1} M_{1} M_{2}) y_{1}^2 y_{2}+ \emph{} \\ \nonumber & \\ \nonumber & \emph{} +(2 C_{2 1 1} M_{1} M_{3}+2 C_{2 1 1} M_{1}^2+12 C_{1 1 1 1} M_{1}^2 M_{3}+4 C_{2 2} M_{1}) y_{1}^2 y_{3}+ \emph{} \\ \nonumber & \\ \nonumber & \emph{} +(2 C_{2 1 1} M_{1} M_{2}+12 C_{1 1 1 1} M_{1} M_{2}^2) y_{1} y_{2}^2+ \emph{} \\ \nonumber & \\ \nonumber & \emph{} +(4 C_{2 1 1} M_{1} M_{2}+24 C_{1 1 1 1} M_{1} M_{2} M_{3}) y_{1} y_{2} y_{3}+ \emph{} \\ \nonumber & \\ \nonumber & \emph{} +(12 C_{1 1 1 1} M_{1} M_{3}^2+3 C_{3 1} M_{1}+6 C_{2 1 1} M_{1} M_{3}) y_{1} y_{3}^2+ \emph{} \\ \nonumber & \\ \nonumber & \emph{} +(4 C_{1 1 1 1} M_{2}^3+C_{3 1} M_{2}+2 C_{2 1 1} M_{2}^2) y_{2}^3+ \emph{} \\ \nonumber & \\ \nonumber & \emph{} +(2 C_{2 1 1} M_{2}^2+2 C_{2 1 1} M_{2} M_{3}+4 C_{2 2} M_{2}+12 C_{1 1 1 1} M_{2}^2 M_{3}) y_{2}^2 y_{3}+ \emph{} \\ \nonumber & \\ \nonumber & \emph{} +(6 C_{2 1 1} M_{2} M_{3}+3 C_{3 1} M_{2}+12 C_{1 1 1 1} M_{2} M_{3}^2) y_{2} y_{3}^2+ \emph{} \\ \nonumber & \\ & \emph{} +(4 C_{2 1 1} M_{3}^2+4 C_{3 1} M_{3}+4 C_{1 1 1 1} M_{3}^3+4 C_{4}+4 C_{2 2} M_{3}) y_{3}^3
\end{align}
\smallskip\\
Their sum equals

\begin{align}
\nonumber \sum\limits_{1 \leq i \leq 3} P^{(i)}_{M_1, M_2, M_3}(y_1, y_2, y_3) \ = \ &
(4 C_{2 2} M_{1}+12 C_{1 1 1 1} M_{1}^3+6 C_{3 1} M_{1}+8 C_{2 1 1} M_{1}^2+4 C_{4}) y_{1}^3+ \emph{} \\ \nonumber & \\ \nonumber & \emph{} +(36 C_{1 1 1 1} M_{1}^2 M_{2}+3 C_{3 1} M_{2}+10 C_{2 1 1} M_{1} M_{2}+2 C_{2 1 1} M_{1}^2+4 C_{2 2} M_{1}) y_{1}^2 y_{2}+ \emph{} \\ \nonumber & \\ \nonumber & \emph{} +(4 C_{2 2} M_{1}+2 C_{2 1 1} M_{1}^2+36 C_{1 1 1 1} M_{1}^2 M_{3}+10 C_{2 1 1} M_{1} M_{3}+3 C_{3 1} M_{3}) y_{1}^2 y_{3}+ \emph{} \\ \nonumber & \\ \nonumber & \emph{} +(3 C_{3 1} M_{1}+36 C_{1 1 1 1} M_{1} M_{2}^2+4 C_{2 2} M_{2}+2 C_{2 1 1} M_{2}^2+10 C_{2 1 1} M_{1} M_{2}) y_{1} y_{2}^2+ \emph{} \\ \nonumber & \\ \nonumber & \emph{} +(4 C_{2 1 1} M_{2} M_{3}+72 C_{1 1 1 1} M_{1} M_{2} M_{3}+4 C_{2 1 1} M_{1} M_{2}+4 C_{2 1 1} M_{1} M_{3}) y_{1} y_{2} y_{3}+ \emph{} \\ \nonumber & \\ \nonumber & \emph{} +(10 C_{2 1 1} M_{1} M_{3}+3 C_{3 1} M_{1}+2 C_{2 1 1} M_{3}^2+4 C_{2 2} M_{3}+36 C_{1 1 1 1} M_{1} M_{3}^2) y_{1} y_{3}^2+ \emph{} \\ \nonumber & \\ \nonumber & \emph{} +(12 C_{1 1 1 1} M_{2}^3+6 C_{3 1} M_{2}+4 C_{2 2} M_{2}+4 C_{4}+8 C_{2 1 1} M_{2}^2) y_{2}^3+ \emph{} \\ \nonumber & \\ \nonumber & \emph{} +(10 C_{2 1 1} M_{2} M_{3}+36 C_{1 1 1 1} M_{2}^2 M_{3}+2 C_{2 1 1} M_{2}^2+4 C_{2 2} M_{2}+3 C_{3 1} M_{3}) y_{2}^2 y_{3}+ \emph{} \\ \nonumber & \\ \nonumber & \emph{} +(36 C_{1 1 1 1} M_{2} M_{3}^2+2 C_{2 1 1} M_{3}^2+3 C_{3 1} M_{2}+4 C_{2 2} M_{3}+10 C_{2 1 1} M_{2} M_{3}) y_{2} y_{3}^2+ \emph{} \\ \nonumber & \\ & \emph{} +(6 C_{3 1} M_{3}+4 C_{2 2} M_{3}+12 C_{1 1 1 1} M_{3}^3+8 C_{2 1 1} M_{3}^2+4 C_{4}) y_{3}^3
\end{align}
\smallskip\\
Using eq. (\ref{2and3levelquantities}), we then calculate the second-level quantities

\begin{align}
\nonumber P^{(1,2)}_{M_1, M_2, M_3}(y_1, y_2, y_3) \ = \ & (4 C_{2 2} M_{1}+3 C_{3 1} M_{1}+4 C_{4}+2 C_{2 1 1} M_{1}^2) y_{1}^2+ \emph{} \\ \nonumber & \\ \nonumber & \emph{} +(4 C_{2 1 1} M_{2} M_{1}+3 C_{3 1} M_{2}+3 C_{3 1} M_{1}+4 C_{4}) y_{1} y_{2}+ \emph{} \\ \nonumber & \\ \nonumber & \emph{} +(4 C_{2 1 1} M_{3} M_{1}+3 C_{3 1} M_{3}) y_{1} y_{3}+ \emph{} \\ \nonumber & \\ \nonumber & \emph{} +(4 C_{2 2} M_{2}+2 C_{2 1 1} M_{2}^2+4 C_{4}+3 C_{3 1} M_{2}) y_{2}^2+ \emph{} \\ \nonumber & \\ \nonumber & \emph{} +(3 C_{3 1} M_{3}+4 C_{2 1 1} M_{2} M_{3}) y_{2} y_{3}+ \emph{} \\ \nonumber & \\ & \emph{} +(2 C_{2 1 1} M_{3}^2+4 C_{2 2} M_{3}) y_{3}^2
\end{align}
\smallskip\\

\begin{align}
\nonumber P^{(1,3)}_{M_1, M_2, M_3}(y_1, y_2, y_3) \ = \ &
(4 C_{2 2} M_{1}+3 C_{3 1} M_{1}+4 C_{4}+2 C_{2 1 1} M_{1}^2) y_{1}^2+ \emph{} \\ \nonumber & \\ \nonumber & \emph{} +(4 C_{2 1 1} M_{2} M_{1}+3 C_{3 1} M_{2}) y_{1} y_{2}+ \emph{} \\ \nonumber & \\ \nonumber & \emph{} +(3 C_{3 1} M_{1}+4 C_{4}+3 C_{3 1} M_{3}+4 C_{2 1 1} M_{3} M_{1}) y_{1} y_{3}+ \emph{} \\ \nonumber & \\ \nonumber & \emph{} +(4 C_{2 2} M_{2}+2 C_{2 1 1} M_{2}^2) y_{2}^2+ \emph{} \\ \nonumber & \\ \nonumber & \emph{} +(3 C_{3 1} M_{2}+4 C_{2 1 1} M_{2} M_{3}) y_{2} y_{3}+ \emph{} \\ \nonumber & \\ & \emph{} +(4 C_{2 2} M_{3}+2 C_{2 1 1} M_{3}^2+3 C_{3 1} M_{3}+4 C_{4}) y_{3}^2
\end{align}

\begin{align}
\nonumber P^{(2,3)}_{M_1, M_2, M_3}(y_1, y_2, y_3) \ = \ &
(4 C_{2 2} M_{1}+2 C_{2 1 1} M_{1}^2) y_{1}^2+ \emph{} \\ \nonumber & \\ \nonumber & \emph{} +(3 C_{3 1} M_{1}+4 C_{2 1 1} M_{2} M_{1}) y_{1} y_{2}+ \emph{} \\ \nonumber & \\ \nonumber & \emph{} +(4 C_{2 1 1} M_{3} M_{1}+3 C_{3 1} M_{1}) y_{1} y_{3}+ \emph{} \\ \nonumber & \\ \nonumber & \emph{} +(4 C_{2 2} M_{2}+2 C_{2 1 1} M_{2}^2+4 C_{4}+3 C_{3 1} M_{2}) y_{2}^2+ \emph{} \\ \nonumber & \\ \nonumber & \emph{} +(3 C_{3 1} M_{3}+4 C_{4}+4 C_{2 1 1} M_{2} M_{3}+3 C_{3 1} M_{2}) y_{2} y_{3}+ \emph{} \\ \nonumber & \\ & \emph{} +(4 C_{2 2} M_{3}+2 C_{2 1 1} M_{3}^2+3 C_{3 1} M_{3}+4 C_{4}) y_{3}^2
\end{align}
\emph{}\vspace{-4ex}

\begin{align}
\nonumber \sum\limits_{1 \leq i < j \leq 3} P^{(i,j)}_{M_1, M_2, M_3}(y_1, y_2, y_3) \ = \ &
(6 C_{3 1} M_{1}+8 C_{4}+12 C_{2 2} M_{1}+6 C_{2 1 1} M_{1}^2) y_{1}^2+ \emph{} \\ \nonumber & \\ \nonumber & \emph{} +(12 C_{2 1 1} M_{2} M_{1}+6 C_{3 1} M_{2}+4 C_{4}+6 C_{3 1} M_{1}) y_{1} y_{2}+ \emph{} \\ \nonumber & \\ \nonumber & \emph{} +(4 C_{4}+6 C_{3 1} M_{1}+6 C_{3 1} M_{3}+12 C_{2 1 1} M_{3} M_{1}) y_{1} y_{3}+ \emph{} \\ \nonumber & \\ \nonumber & \emph{} +(12 C_{2 2} M_{2}+8 C_{4}+6 C_{3 1} M_{2}+6 C_{2 1 1} M_{2}^2) y_{2}^2+ \emph{} \\ \nonumber & \\ \nonumber & \emph{} +(12 C_{2 1 1} M_{2} M_{3}+4 C_{4}+6 C_{3 1} M_{3}+6 C_{3 1} M_{2}) y_{2} y_{3}+ \emph{} \\ \nonumber & \\ & \emph{} +(6 C_{3 1} M_{3}+8 C_{4}+6 C_{2 1 1} M_{3}^2+12 C_{2 2} M_{3}) y_{3}^2
\end{align}
\smallskip\\
Finally, we calculate the single third-level quantity:

\begin{align}
P^{(1,2,3)}_{M_1, M_2, M_3}(y_1, y_2, y_3) = (3 C_{3 1} M_{1}+4 C_{4}) y_{1}+(3 C_{3 1} M_{2}+4 C_{4}) y_{2}+(4 C_{4}+3 C_{3 1} M_{3}) y_{3}
\end{align}
\smallskip\\
According to (\ref{MainFormula}), the irreducible factors are equal to

\begin{align}
d_{M_1 M_2 M_3}\big(C\big) = R\Big\{ \sum\limits_{1 \leq i \leq 3} P^{(i)}_{M_1, M_2, M_3}, \sum\limits_{1 \leq i < j \leq 3} P^{(i,j)}_{M_1, M_2, M_3}, P^{(1,2,3)}_{M_1, M_2, M_3} \Big\}
\end{align}
\smallskip\\
In the right hand side we have a resultant of three homogeneous polynomials of degrees $3,2,1$ in variables $y_1, y_2, y_3$. Such a resultant can be computed by solving the last linear polynomial equation, substituting into the first two and then taking a Sylvester resultant. After a tedious computer calculation, one obtains

\begin{align*}
d_{M_1 M_2 M_3}\big( C \big) =
\end{align*}

\begin{center}
$=65536 C_{4}^9 (C_{1111} n^3+C_{211} n^2+C_{22} n+C_{31} n+C_{4})^2 + 243 (C_{31}^4 C_{22} n-256 C_{4}^3 C_{22} C_{1111}+9 C_{4} C_{31}^4-48 C_{211} C_{4}^2 C_{31}^2+64 C_{4}^3 C_{211}^2+16 C_{4} C_{31}^3 C_{22}) (3 C_{31}^4 C_{22}^2 n^2+18 C_{4} C_{31}^4 C_{22} n+48 C_{4} C_{31}^3 C_{22}^2 n-80 C_{211} C_{4}^2 C_{31}^2 C_{22} n-256 C_{4}^3 C_{22}^2 C_{1111} n+64 C_{4}^3 C_{211}^2 C_{22} n-9 C_{4}^2 C_{31}^4+48 C_{4}^2 C_{31}^3 C_{22}+128 C_{4}^2 C_{31}^2 C_{22}^2+48 C_{211} C_{4}^3 C_{31}^2-256 C_{211} C_{4}^3 C_{31} C_{22}+768 C_{4}^4 C_{22} C_{1111}-64 C_{4}^4 C_{211}^2) \sigma_{3}^2 \sigma_{2}+72 C_{4}^3 (336 C_{4} C_{31}^4 C_{22} C_{211}^2 n^4-512 C_{4}^2 C_{31}^2 C_{211}^3 C_{22} n^4-54 C_{31}^6 C_{22} C_{211} n^4+2048 C_{4}^2 C_{31}^2 C_{22}^2 C_{1111} C_{211} n^4+1344 C_{4} C_{31}^4 C_{22}^2 C_{211} n^3+6144 C_{4}^2 C_{31}^3 C_{22}^2 C_{1111} n^3-21504 C_{4}^3 C_{31}^2 C_{22} C_{1111} C_{211} n^3+768 C_{4}^2 C_{31}^3 C_{211}^2 C_{22} n^3+16384 C_{4}^4 C_{22} C_{1111} C_{211}^2 n^3+5184 C_{4}^2 C_{31}^4 C_{22} C_{1111} n^3+3584 C_{4}^3 C_{31}^2 C_{211}^3 n^3-162 C_{31}^7 C_{22} n^3+4096 C_{4}^2 C_{31}^2 C_{22}^3 C_{1111} n^3+720 C_{4} C_{31}^5 C_{22} C_{211} n^3-2784 C_{4}^2 C_{31}^4 C_{211}^2 n^3-1024 C_{4}^4 C_{211}^4 n^3+16384 C_{4}^3 C_{31} C_{22}^2 C_{211} C_{1111} n^3-1024 C_{4}^2 C_{31}^2 C_{22}^2 C_{211}^2 n^3-81 C_{31}^8 n^3-4096 C_{4}^3 C_{31} C_{211}^3 C_{22} n^3-108 C_{31}^6 C_{22}^2 n^3-49152 C_{4}^4 C_{1111}^2 C_{22}^2 n^3+810 C_{4} C_{31}^6 C_{211} n^3+18432 C_{4}^3 C_{31}^2 C_{22}^2 C_{1111} n^2-24576 C_{4}^5 C_{211}^2 C_{1111} n^2+13824 C_{4}^4 C_{31}^2 C_{1111} C_{211} n^2+16384 C_{4}^4 C_{31} C_{211}^3 n^2+2736 C_{4}^2 C_{31}^5 C_{211} n^2+3552 C_{4}^2 C_{31}^4 C_{22} C_{211} n^2-73728 C_{4}^4 C_{31} C_{22} C_{1111} C_{211} n^2-1728 C_{31}^4 C_{1111} C_{4}^3 n^2+32768 C_{4}^3 C_{31} C_{22}^3 C_{1111} n^2+432 C_{4} C_{31}^6 C_{22} n^2-12288 C_{4}^3 C_{31}^3 C_{211}^2 n^2+18432 C_{4}^3 C_{31}^3 C_{22} C_{1111} n^2+1344 C_{4} C_{31}^4 C_{22}^3 n^2-162 C_{4} C_{31}^7 n^2+11264 C_{4}^2 C_{31}^3 C_{22}^2 C_{211} n^2-8192 C_{4}^3 C_{31} C_{22}^2 C_{211}^2 n^2+147456 C_{4}^5 C_{22} C_{1111}^2 n^2-13824 C_{4}^3 C_{31}^2 C_{22} C_{211}^2 n^2+1440 C_{4} C_{31}^5 C_{22}^2 n^2-8192 C_{4}^4 C_{22}^2 C_{211}^2 n-11232 C_{4}^3 C_{31}^4 C_{211} n+10656 C_{4}^2 C_{31}^5 C_{22} n+27648 C_{4}^4 C_{31}^2 C_{22} C_{1111} n-20480 C_{4}^4 C_{31} C_{22} C_{211}^2 n+1836 C_{4}^2 C_{31}^6 n+13312 C_{4}^2 C_{31}^3 C_{22}^3 n-49152 C_{4}^4 C_{31} C_{22}^2 C_{1111} n+19200 C_{4}^2 C_{31}^4 C_{22}^2 n-36864 C_{4}^5 C_{31} C_{1111} C_{211} n+32768 C_{4}^4 C_{22}^3 C_{1111} n-29184 C_{4}^3 C_{31}^3 C_{22} C_{211} n+49152 C_{4}^5 C_{22} C_{1111} C_{211} n+13824 C_{4}^4 C_{31}^2 C_{211}^2 n+3072 C_{4}^3 C_{31}^2 C_{22}^2 C_{211} n+13824 C_{31}^3 C_{1111} C_{4}^4 n+8192 C_{4}^5 C_{211}^3 n-73728 C_{4}^6 C_{211} C_{1111}+36864 C_{4}^5 C_{31} C_{211}^2+28672 C_{4}^3 C_{31}^2 C_{22}^3-23040 C_{4}^4 C_{31}^3 C_{211}-57344 C_{4}^4 C_{31} C_{22}^2 C_{211}+73728 C_{4}^5 C_{31} C_{22} C_{1111}+39936 C_{31}^3 C_{22}^2 C_{4}^3+20160 C_{31}^4 C_{22} C_{4}^3-76800 C_{4}^4 C_{31}^2 C_{211} C_{22}+27648 C_{31}^2 C_{1111} C_{4}^5+28672 C_{4}^5 C_{211}^2 C_{22}+3456 C_{4}^3 C_{31}^5) \sigma_{3} \sigma_{2}-432 C_{4}^3 (32 C_{31}^4 C_{22}^3 C_{1111} n^3-8 C_{31}^4 C_{22}^2 C_{211}^2 n^3+1536 C_{4}^2 C_{31}^2 C_{22}^2 C_{1111} C_{211} n^2+96 C_{4} C_{31}^4 C_{22}^2 C_{1111} n^2+4096 C_{4}^3 C_{22}^3 C_{1111}^2 n^2-36 C_{31}^6 C_{22} C_{211} n^2-48 C_{31}^5 C_{22}^2 C_{211} n^2-128 C_{4} C_{31}^3 C_{22}^2 C_{211}^2 n^2+208 C_{4} C_{31}^4 C_{22} C_{211}^2 n^2-384 C_{4}^2 C_{31}^2 C_{211}^3 C_{22} n^2-2048 C_{4}^3 C_{22}^2 C_{1111} C_{211}^2 n^2+512 C_{4} C_{31}^3 C_{22}^3 C_{1111} n^2+256 C_{4}^3 C_{22} C_{211}^4 n^2+1152 C_{4}^3 C_{31}^2 C_{211}^3 n+252 C_{4} C_{31}^6 C_{211} n-64 C_{31}^5 C_{22}^3 n-27 C_{31}^8 n-512 C_{4}^4 C_{211}^4 n+64 C_{4} C_{31}^4 C_{22}^2 C_{211} n-144 C_{31}^6 C_{22}^2 n+2016 C_{4}^2 C_{31}^4 C_{22} C_{1111} n+384 C_{4} C_{31}^5 C_{22} C_{211} n+6144 C_{4}^2 C_{31}^3 C_{22}^2 C_{1111} n+6144 C_{4}^2 C_{31}^2 C_{22}^3 C_{1111} n+8192 C_{4}^4 C_{22} C_{1111} C_{211}^2 n-128 C_{4}^2 C_{31}^3 C_{211}^2 C_{22} n-24576 C_{4}^4 C_{1111}^2 C_{22}^2 n-1536 C_{4}^2 C_{31}^2 C_{22}^2 C_{211}^2 n-9216 C_{4}^3 C_{31}^2 C_{22} C_{1111} C_{211} n-108 C_{31}^7 C_{22} n-840 C_{4}^2 C_{31}^4 C_{211}^2 n+3072 C_{4}^4 C_{31} C_{211}^3+2112 C_{4}^2 C_{31}^4 C_{22} C_{211}-4096 C_{4}^3 C_{31} C_{22}^2 C_{211}^2-3072 C_{4}^3 C_{31}^2 C_{22} C_{211}^2+1024 C_{4}^2 C_{31}^3 C_{22}^2 C_{211}+4608 C_{4}^3 C_{31}^3 C_{22} C_{1111}+16384 C_{4}^3 C_{31} C_{22}^3 C_{1111}-864 C_{31}^4 C_{1111} C_{4}^3+4608 C_{4}^4 C_{31}^2 C_{1111} C_{211}+36864 C_{4}^5 C_{22} C_{1111}^2-24576 C_{4}^4 C_{31} C_{22} C_{1111} C_{211}+18432 C_{31}^2 C_{22}^2 C_{4}^3 C_{1111}+1008 C_{4}^2 C_{31}^5 C_{211}-3072 C_{4}^3 C_{31}^3 C_{211}^2-432 C_{31}^6 C_{22} C_{4}-576 C_{31}^5 C_{22}^2 C_{4}-16384 C_{4}^4 C_{22}^2 C_{211} C_{1111}-256 C_{4} C_{31}^4 C_{22}^3-6144 C_{4}^5 C_{211}^2 C_{1111}+4096 C_{4}^4 C_{22} C_{211}^3-108 C_{4} C_{31}^7) \sigma_{3} \sigma_{2}^2+256 C_{4}^6 (C_{1111} n^3+C_{211} n^2+C_{22} n+C_{31} n+C_{4}) (768 C_{4}^2 C_{1111} C_{22} n^2-27 C_{31}^4 n^2-128 C_{4}^2 C_{211}^2 n^2+144 C_{211} C_{31}^2 C_{4} n^2-2304 C_{4}^3 C_{1111} n+384 C_{211} C_{31} C_{4}^2 n+288 C_{4} C_{31}^2 C_{22} n+256 C_{211} C_{4}^2 C_{22} n+256 C_{22}^2 C_{4}^2+576 C_{4}^2 C_{31}^2-1536 C_{211} C_{4}^3+1536 C_{4}^2 C_{31} C_{22}) \sigma_{2}+C_{4}^3 (27648 C_{4}^2 C_{31}^4 C_{211}^2 n^5+729 C_{31}^8 n^5+147456 C_{4}^3 C_{31}^2 C_{22} C_{1111} C_{211} n^5-41472 C_{4}^2 C_{31}^4 C_{22} C_{1111} n^5-32768 C_{4}^3 C_{31}^2 C_{211}^3 n^5-7776 C_{4} C_{31}^6 C_{211} n^5+69120 C_{4}^2 C_{31}^4 C_{22} C_{211} n^4+294912 C_{4}^5 C_{211}^2 C_{1111} n^4-20736 C_{4}^2 C_{31}^5 C_{211} n^4+147456 C_{4}^3 C_{31}^3 C_{211}^2 n^4-1769472 C_{4}^5 C_{22} C_{1111}^2 n^4+294912 C_{4}^3 C_{31}^2 C_{22}^2 C_{1111} n^4-49152 C_{4}^3 C_{31}^2 C_{22} C_{211}^2 n^4-110592 C_{4}^4 C_{31}^2 C_{1111} C_{211} n^4+6912 C_{4}^3 C_{31}^4 C_{1111} n^4-15552 C_{4} C_{31}^6 C_{22} n^4+1179648 C_{4}^4 C_{31} C_{22} C_{1111} C_{211} n^4-221184 C_{4}^3 C_{31}^3 C_{22} C_{1111} n^4-262144 C_{4}^4 C_{31} C_{211}^3 n^4+69120 C_{4}^2 C_{31}^4 C_{22}^2 n^3+221184 C_{4}^4 C_{31}^2 C_{22} C_{1111} n^3-405504 C_{4}^4 C_{31}^2 C_{211}^2 n^3+49152 C_{4}^3 C_{31}^2 C_{22}^2 C_{211} n^3-229376 C_{4}^5 C_{211}^3 n^3+241920 C_{4}^3 C_{31}^4 C_{211} n^3-31104 C_{4}^2 C_{31}^6 n^3-82944 C_{4}^2 C_{31}^5 C_{22} n^3+884736 C_{4}^5 C_{31} C_{1111} C_{211} n^3-393216 C_{4}^4 C_{31} C_{22} C_{211}^2 n^3+516096 C_{4}^3 C_{31}^3 C_{22} C_{211} n^3-221184 C_{4}^4 C_{31}^3 C_{1111} n^3+2359296 C_{4}^4 C_{31} C_{22}^2 C_{1111} n^3+1769472 C_{4}^5 C_{31} C_{22} C_{1111} n^2-1179648 C_{4}^5 C_{31} C_{211}^2 n^2+214272 C_{4}^3 C_{31}^4 C_{22} n^2+258048 C_{4}^4 C_{31}^2 C_{211} C_{22} n^2-6912 C_{4}^3 C_{31}^5 n^2+32768 C_{4}^3 C_{31}^2 C_{22}^3 n^2-1081344 C_{4}^5 C_{211}^2 C_{22} n^2+663552 C_{4}^3 C_{31}^3 C_{22}^2 n^2+393216 C_{4}^4 C_{31} C_{22}^2 C_{211} n^2+331776 C_{4}^4 C_{31}^3 C_{211} n^2+2359296 C_{4}^5 C_{22}^2 C_{1111} n^2+1769472 C_{4}^6 C_{211} C_{1111} n^2-2064384 C_{4}^5 C_{31} C_{22} C_{211} n+262144 C_{4}^4 C_{31} C_{22}^3 n+1769472 C_{4}^6 C_{31} C_{1111} n-393216 C_{4}^5 C_{22}^2 C_{211} n+1105920 C_{4}^4 C_{31}^3 C_{22} n-294912 C_{4}^6 C_{211}^2 n+1695744 C_{4}^4 C_{31}^2 C_{22}^2 n+214272 C_{4}^4 C_{31}^4 n+3538944 C_{4}^6 C_{22} C_{1111} n-774144 C_{4}^5 C_{31}^2 C_{211} n-884736 C_{4}^6 C_{31} C_{211}+589824 C_{4}^5 C_{31} C_{22}^2-589824 C_{211} C_{4}^6 C_{22}+1769472 C_{1111} C_{4}^7-65536 C_{4}^5 C_{22}^3+221184 C_{4}^5 C_{31}^3+663552 C_{31}^2 C_{4}^5 C_{22}) \sigma_{3}+1024 C_{4}^6 (-C_{31}^2 C_{211}^3 n^4+192 C_{4} C_{1111}^2 C_{22}^2 n^4+4 C_{4} C_{211}^4 n^4-64 C_{4} C_{22} C_{1111} C_{211}^2 n^4+36 C_{31}^2 C_{22} C_{1111} C_{211} n^4-32 C_{4} C_{22} C_{211}^3 n^3-96 C_{4} C_{31} C_{22} C_{1111} C_{211} n^3+192 C_{4}^2 C_{211}^2 C_{1111} n^3-1152 C_{4}^2 C_{22} C_{1111}^2 n^3+54 C_{31}^4 C_{1111} n^3+40 C_{4} C_{31} C_{211}^3 n^3+72 C_{31}^2 C_{22}^2 C_{1111} n^3-9 C_{31}^3 C_{211}^2 n^3+108 C_{31}^3 C_{22} C_{1111} n^3-252 C_{4} C_{31}^2 C_{1111} C_{211} n^3+30 C_{31}^2 C_{22} C_{211}^2 n^3+128 C_{4} C_{22}^2 C_{211} C_{1111} n^3+128 C_{4} C_{22}^3 C_{1111} n^2-126 C_{4} C_{31}^2 C_{211}^2 n^2+108 C_{31}^3 C_{22} C_{211} n^2-32 C_{4} C_{22}^2 C_{211}^2 n^2-864 C_{4}^2 C_{31} C_{1111} C_{211} n^2+96 C_{31}^2 C_{22}^2 C_{211} n^2+192 C_{4} C_{31} C_{22}^2 C_{1111} n^2+108 C_{4} C_{31}^3 C_{1111} n^2+272 C_{4}^2 C_{211}^3 n^2+80 C_{4} C_{31} C_{22} C_{211}^2 n^2-1728 C_{4}^2 C_{22} C_{1111} C_{211} n^2+1728 C_{4}^3 C_{1111}^2 n^2+27 C_{31}^4 C_{211} n^2+32 C_{4}^2 C_{211}^2 C_{22} n+27 C_{31}^5 n-36 C_{4} C_{31}^3 C_{211} n+64 C_{31}^2 C_{22}^3 n+192 C_{4} C_{31}^2 C_{211} C_{22} n-96 C_{4}^2 C_{31} C_{211}^2 n+1728 C_{4}^3 C_{211} C_{1111} n-2880 C_{4}^2 C_{31} C_{22} C_{1111} n-648 C_{4}^2 C_{31}^2 C_{1111} n+108 C_{31}^4 C_{22} n+144 C_{31}^3 C_{22}^2 n+320 C_{4} C_{31} C_{22}^2 C_{211} n-1536 C_{4}^2 C_{22}^2 C_{1111} n+256 C_{4} C_{31} C_{22}^3-1152 C_{4}^3 C_{22} C_{1111}+432 C_{4} C_{31}^3 C_{22}+108 C_{4} C_{31}^4-960 C_{4}^2 C_{31} C_{22} C_{211}-256 C_{4}^2 C_{22}^2 C_{211}+576 C_{4} C_{31}^2 C_{22}^2-576 C_{211} C_{4}^2 C_{31}^2+768 C_{4}^3 C_{211}^2) \sigma_{2}^2+4096 C_{4}^6 (C_{22} C_{211}^4 n^3-8 C_{22}^2 C_{1111} C_{211}^2 n^3+16 C_{22}^3 C_{1111}^2 n^3-3 C_{4} C_{211}^4 n^2-144 C_{4} C_{1111}^2 C_{22}^2 n^2+48 C_{4} C_{22} C_{1111} C_{211}^2 n^2-36 C_{31}^2 C_{22} C_{1111} C_{211} n^2-48 C_{31} C_{22}^2 C_{211} C_{1111} n^2+C_{31}^2 C_{211}^3 n^2+12 C_{31} C_{211}^3 C_{22} n^2+9 C_{31}^3 C_{211}^2 n-72 C_{4}^2 C_{211}^2 C_{1111} n+48 C_{4} C_{22} C_{211}^3 n+108 C_{4} C_{31}^2 C_{1111} C_{211} n-144 C_{31}^2 C_{22}^2 C_{1111} n+24 C_{31}^2 C_{22} C_{211}^2 n+16 C_{31} C_{22}^2 C_{211}^2 n+432 C_{4}^2 C_{22} C_{1111}^2 n-27 C_{31}^4 C_{1111} n-28 C_{4} C_{31} C_{211}^3 n-64 C_{31} C_{22}^3 C_{1111} n-108 C_{31}^3 C_{22} C_{1111} n-192 C_{4} C_{22}^2 C_{211} C_{1111} n+96 C_{4} C_{31} C_{22} C_{211}^2-432 C_{4} C_{31}^2 C_{22} C_{1111}-108 C_{4} C_{31}^3 C_{1111}-576 C_{4} C_{31} C_{22}^2 C_{1111}+64 C_{4} C_{22}^2 C_{211}^2+432 C_{4}^2 C_{31} C_{1111} C_{211}+36 C_{4} C_{31}^2 C_{211}^2-432 C_{4}^3 C_{1111}^2+576 C_{4}^2 C_{22} C_{1111} C_{211}-256 C_{4} C_{22}^3 C_{1111}-128 C_{4}^2 C_{211}^3) \sigma_{2}^3-729 C_{22} (C_{31}^4 C_{22} n-256 C_{4}^3 C_{22} C_{1111}+9 C_{4} C_{31}^4-48 C_{211} C_{4}^2 C_{31}^2+64 C_{4}^3 C_{211}^2+16 C_{4} C_{31}^3 C_{22})^2 \sigma_{3}^3+27 C_{4} (-1024 C_{4}^2 C_{31}^4 C_{22}^2 C_{211}^2 n^4+4096 C_{4}^2 C_{31}^4 C_{22}^3 C_{1111} n^4+576 C_{4} C_{31}^6 C_{211} C_{22}^2 n^4-81 C_{31}^8 C_{22}^2 n^4+18432 C_{4}^3 C_{31}^4 C_{22}^2 C_{1111} n^3-11136 C_{4}^3 C_{31}^4 C_{22} C_{211}^2 n^3-16384 C_{4}^3 C_{31}^3 C_{22}^2 C_{211}^2 n^3-49152 C_{4}^4 C_{31}^2 C_{22}^2 C_{1111} C_{211} n^3+7680 C_{4}^2 C_{31}^5 C_{22}^2 C_{211} n^3+1152 C_{4} C_{31}^6 C_{22}^3 n^3-864 C_{4} C_{31}^7 C_{22}^2 n^3+3312 C_{4}^2 C_{31}^6 C_{22} C_{211} n^3-324 C_{4} C_{31}^8 C_{22} n^3+65536 C_{4}^3 C_{31}^3 C_{22}^3 C_{1111} n^3+12288 C_{4}^4 C_{31}^2 C_{211}^3 C_{22} n^3+1536 C_{4}^3 C_{31}^4 C_{22}^2 C_{211} n^2-73728 C_{4}^4 C_{31}^2 C_{22}^2 C_{211}^2 n^2+189 C_{4}^2 C_{31}^8 n^2-13824 C_{4}^4 C_{31}^4 C_{22} C_{1111} n^2-393216 C_{4}^5 C_{31} C_{22}^2 C_{211} C_{1111} n^2+9216 C_{4}^4 C_{31}^4 C_{211}^2 n^2-17408 C_{4}^5 C_{31}^2 C_{211}^3 n^2+110592 C_{4}^5 C_{31}^2 C_{22} C_{1111} C_{211} n^2-67584 C_{4}^4 C_{31}^3 C_{211}^2 C_{22} n^2+294912 C_{4}^4 C_{31}^2 C_{22}^3 C_{1111} n^2-432 C_{4}^2 C_{31}^7 C_{22} n^2+147456 C_{4}^4 C_{31}^3 C_{22}^2 C_{1111} n^2+12672 C_{4}^3 C_{31}^5 C_{22} C_{211} n^2+19456 C_{4}^2 C_{31}^5 C_{22}^3 n^2+98304 C_{4}^5 C_{31} C_{211}^3 C_{22} n^2-196608 C_{4}^6 C_{22} C_{1111} C_{211}^2 n^2+12288 C_{4}^6 C_{211}^4 n^2+6048 C_{4}^2 C_{31}^6 C_{22}^2 n^2-2160 C_{4}^3 C_{31}^6 C_{211} n^2+589824 C_{4}^6 C_{1111}^2 C_{22}^2 n^2-294912 C_{4}^6 C_{31} C_{22} C_{1111} C_{211} n+66816 C_{4}^3 C_{31}^5 C_{22}^2 n+98304 C_{4}^6 C_{22} C_{211}^3 n-432 C_{4}^3 C_{31}^7 n-9216 C_{4}^5 C_{31}^3 C_{211}^2 n+98816 C_{4}^3 C_{31}^4 C_{22}^3 n-221184 C_{4}^4 C_{31}^3 C_{22}^2 C_{211} n-65536 C_{4}^5 C_{31} C_{22}^2 C_{211}^2 n+516096 C_{4}^5 C_{31}^2 C_{22}^2 C_{1111} n+43008 C_{4}^5 C_{31}^2 C_{22} C_{211}^2 n+8640 C_{4}^3 C_{31}^6 C_{22} n-393216 C_{4}^6 C_{22}^2 C_{211} C_{1111} n+110592 C_{4}^5 C_{31}^3 C_{22} C_{1111} n-52992 C_{4}^4 C_{31}^4 C_{22} C_{211} n+262144 C_{4}^5 C_{31} C_{22}^3 C_{1111} n+8192 C_{4}^6 C_{31} C_{211}^3 n+3456 C_{4}^4 C_{31}^5 C_{211} n+20736 C_{4}^4 C_{31}^5 C_{22}+139264 C_{4}^4 C_{31}^3 C_{22}^3+6912 C_{4}^5 C_{31}^4 C_{211}+245760 C_{4}^6 C_{31} C_{22} C_{211}^2+115200 C_{4}^4 C_{31}^4 C_{22}^2-18432 C_{4}^6 C_{31}^2 C_{211}^2-147456 C_{4}^5 C_{31}^3 C_{22} C_{211}+16384 C_{4}^7 C_{211}^3+589824 C_{4}^6 C_{31} C_{22}^2 C_{1111}-417792 C_{4}^5 C_{31}^2 C_{22}^2 C_{211}+221184 C_{4}^6 C_{31}^2 C_{22} C_{1111}-589824 C_{4}^7 C_{22} C_{1111} C_{211}-131072 C_{4}^6 C_{22}^3 C_{1111}+32768 C_{4}^6 C_{22}^2 C_{211}^2-864 C_{4}^4 C_{31}^6) \sigma_{3}^2$
\end{center}
where $n = \sigma_1 = M_1 + M_2 + M_3, \ \sigma_2 = M_1 M_2 + M_1 M_3 + M_2 M_3$ and $\sigma_3 = M_1 M_2 M_3$. Similarly, using the formulas from section 2.4, one can obtain the expressions when one of $M_i$ vanishes

\begin{align*}
d_{M_1 M_2 0}\big( C\big) = d_{M_1 0 M_2}\big( C\big) = d_{0 M_1 M_2}\big( C\big) =
\end{align*}

\begin{center}
$
= 4096 C_{4}^3 (C_{1 1 1 1} n^3+C_{2 1 1} n^2+C_{2 2} n+C_{3 1} n+C_{4})^2+16 (C_{1 1 1 1} n^3+C_{2 1 1} n^2+C_{2 2} n+C_{3 1} n+C_{4}) (768 C_{1 1 1 1} C_{4}^2 C_{2 2} n^2+144 C_{3 1}^2 C_{4} C_{2 1 1} n^2-128 C_{4}^2 C_{2 1 1}^2 n^2-27 C_{3 1}^4 n^2-2304 C_{4}^3 C_{1 1 1 1} n+384 C_{4}^2 C_{2 1 1} C_{3 1} n+256 C_{4}^2 C_{2 1 1} C_{2 2} n+288 C_{4} C_{3 1}^2 C_{2 2} n+576 C_{3 1}^2 C_{4}^2-1536 C_{4}^3 C_{2 1 1}+256 C_{4}^2 C_{2 2}^2+1536 C_{2 2} C_{3 1} C_{4}^2) M_{1} M_{2}+(12288 C_{4} C_{1 1 1 1}^2 C_{2 2}^2 n^4-4096 C_{2 1 1}^2 C_{1 1 1 1} C_{4} C_{2 2} n^4+256 C_{2 1 1}^4 C_{4} n^4-64 C_{3 1}^2 C_{2 1 1}^3 n^4+2304 C_{3 1}^2 C_{1 1 1 1} C_{2 2} C_{2 1 1} n^4+12288 C_{2 1 1}^2 C_{4}^2 C_{1 1 1 1} n^3-6144 C_{2 1 1} C_{1 1 1 1} C_{4} C_{3 1} C_{2 2} n^3-576 C_{2 1 1}^2 C_{3 1}^3 n^3-16128 C_{4} C_{1 1 1 1} C_{3 1}^2 C_{2 1 1} n^3-73728 C_{4}^2 C_{1 1 1 1}^2 C_{2 2} n^3+3456 C_{3 1}^4 C_{1 1 1 1} n^3+8192 C_{4} C_{2 2}^2 C_{1 1 1 1} C_{2 1 1} n^3+4608 C_{3 1}^2 C_{1 1 1 1} C_{2 2}^2 n^3+6912 C_{3 1}^3 C_{1 1 1 1} C_{2 2} n^3+2560 C_{2 1 1}^3 C_{4} C_{3 1} n^3+1920 C_{3 1}^2 C_{2 1 1}^2 C_{2 2} n^3-2048 C_{2 1 1}^3 C_{2 2} C_{4} n^3+6912 C_{2 2} C_{2 1 1} C_{3 1}^3 n^2+110592 C_{1 1 1 1}^2 C_{4}^3 n^2-8064 C_{3 1}^2 C_{2 1 1}^2 C_{4} n^2+6144 C_{2 2}^2 C_{2 1 1} C_{3 1}^2 n^2+5120 C_{2 2} C_{2 1 1}^2 C_{4} C_{3 1} n^2-2048 C_{2 2}^2 C_{2 1 1}^2 C_{4} n^2+17408 C_{4}^2 C_{2 1 1}^3 n^2+6912 C_{4} C_{1 1 1 1} C_{3 1}^3 n^2-55296 C_{4}^2 C_{2 1 1} C_{1 1 1 1} C_{3 1} n^2+12288 C_{4} C_{1 1 1 1} C_{2 2}^2 C_{3 1} n^2-110592 C_{2 2} C_{2 1 1} C_{1 1 1 1} C_{4}^2 n^2+8192 C_{4} C_{2 2}^3 C_{1 1 1 1} n^2+1728 C_{3 1}^4 C_{2 1 1} n^2+6912 C_{3 1}^4 C_{2 2} n+110592 C_{4}^3 C_{1 1 1 1} C_{2 1 1} n-184320 C_{4}^2 C_{3 1} C_{1 1 1 1} C_{2 2} n-98304 C_{4}^2 C_{2 2}^2 C_{1 1 1 1} n+4096 C_{2 2}^3 C_{3 1}^2 n-2304 C_{3 1}^3 C_{4} C_{2 1 1} n-41472 C_{3 1}^2 C_{1 1 1 1} C_{4}^2 n-6144 C_{2 1 1}^2 C_{3 1} C_{4}^2 n+12288 C_{4} C_{2 2} C_{3 1}^2 C_{2 1 1} n+2048 C_{2 1 1}^2 C_{4}^2 C_{2 2} n+9216 C_{3 1}^3 C_{2 2}^2 n+20480 C_{2 1 1} C_{3 1} C_{4} C_{2 2}^2 n+1728 C_{3 1}^5 n+49152 C_{4}^3 C_{2 1 1}^2+27648 C_{4} C_{2 2} C_{3 1}^3-36864 C_{3 1}^2 C_{4}^2 C_{2 1 1}+16384 C_{4} C_{2 2}^3 C_{3 1}+6912 C_{3 1}^4 C_{4}-16384 C_{2 1 1} C_{4}^2 C_{2 2}^2-73728 C_{4}^3 C_{2 2} C_{1 1 1 1}+36864 C_{2 2}^2 C_{3 1}^2 C_{4}-61440 C_{4}^2 C_{2 1 1} C_{3 1} C_{2 2}) M_{1}^2 M_{2}^2+(256 C_{2 2} C_{2 1 1}^4 n^3-2048 C_{2 2}^2 C_{2 1 1}^2 C_{1 1 1 1} n^3+4096 C_{2 2}^3 C_{1 1 1 1}^2 n^3-12288 C_{2 1 1} C_{3 1} C_{1 1 1 1} C_{2 2}^2 n^2+3072 C_{2 1 1}^3 C_{3 1} C_{2 2} n^2-768 C_{2 1 1}^4 C_{4} n^2+12288 C_{2 1 1}^2 C_{1 1 1 1} C_{4} C_{2 2} n^2-9216 C_{3 1}^2 C_{1 1 1 1} C_{2 2} C_{2 1 1} n^2-36864 C_{4} C_{1 1 1 1}^2 C_{2 2}^2 n^2+256 C_{3 1}^2 C_{2 1 1}^3 n^2-49152 C_{4} C_{2 2}^2 C_{1 1 1 1} C_{2 1 1} n-27648 C_{3 1}^3 C_{1 1 1 1} C_{2 2} n+2304 C_{2 1 1}^2 C_{3 1}^3 n+6144 C_{3 1}^2 C_{2 1 1}^2 C_{2 2} n-18432 C_{1 1 1 1} C_{4}^2 C_{2 1 1}^2 n+110592 C_{1 1 1 1}^2 C_{4}^2 C_{2 2} n-36864 C_{3 1}^2 C_{1 1 1 1} C_{2 2}^2 n+4096 C_{2 1 1}^2 C_{2 2}^2 C_{3 1} n-16384 C_{2 2}^3 C_{3 1} C_{1 1 1 1} n+12288 C_{2 1 1}^3 C_{2 2} C_{4} n-7168 C_{2 1 1}^3 C_{4} C_{3 1} n-6912 C_{1 1 1 1} C_{3 1}^4 n+27648 C_{1 1 1 1} C_{3 1}^2 C_{4} C_{2 1 1} n-110592 C_{1 1 1 1}^2 C_{4}^3+16384 C_{2 2}^2 C_{2 1 1}^2 C_{4}-32768 C_{4}^2 C_{2 1 1}^3+24576 C_{2 2} C_{2 1 1}^2 C_{4} C_{3 1}-110592 C_{2 2} C_{3 1}^2 C_{1 1 1 1} C_{4}-147456 C_{4} C_{1 1 1 1} C_{2 2}^2 C_{3 1}-27648 C_{4} C_{1 1 1 1} C_{3 1}^3+9216 C_{3 1}^2 C_{2 1 1}^2 C_{4}-65536 C_{4} C_{2 2}^3 C_{1 1 1 1}+147456 C_{2 2} C_{2 1 1} C_{1 1 1 1} C_{4}^2+110592 C_{4}^2 C_{2 1 1} C_{1 1 1 1} C_{3 1}) M_{1}^3 M_{2}^3
$\smallskip\\
\end{center}
and when a pair of $M_i$ vanishes:

\begin{align*}
d_{n 0 0}\big( C\big) = d_{0 n 0}\big( C\big) = d_{0 0 n}\big( C\big) = 4(C_{1 1 1 1} n^3+C_{2 1 1} n^2+C_{2 2} n+C_{3 1} n+C_{4})
\end{align*}
\smallskip\\
So we have calculated all the necessary resultants for the case $r = 4$. Due to their large size, the formula (\ref{4Main}) can not be written as simply as the formula (\ref{DegreeThree}) or (\ref{DegreeThree2}). It can be rather understood as an explicit computer algorithm. Expressions for particular $n$ are straightforward to obtain: say, for $n = 4$ we get
\smallskip\\

\begin{center}
${\cal D}_{4|4}\big(C_{4}, C_{31}, C_{22}, C_{211}, C_{1111}\big) = 2^{-106} C_4^{20} (64 C_{1 1 1 1}+16 C_{2 1 1}+4 C_{3 1}+4 C_{2 2}+C_{4}) (4 C_{22}+C_{4})^3 (2 C_{22}+C_{4})^6  (64 C_{1111} C_{22}+16 C_{1111} C_{4}-16 C_{211}^2-24 C_{211} C_{31}-8 C_{211} C_{4}-9 C_{31}^2-8 C_{31} C_{22}-8 C_{31} C_{4}-4 C_{22} C_{4}-2 C_{4}^2)^6 (256 C_{1111} C_{22} C_{4}^3+128 C_{1111} C_{4}^4-64 C_{211}^2 C_{4}^3+144 C_{211} C_{31}^2 C_{4}^2+192 C_{211} C_{31} C_{4}^3+128 C_{211} C_{4}^4+324 C_{31}^4 C_{22}+81 C_{31}^4 C_{4}+1296 C_{31}^3 C_{22} C_{4}+432 C_{31}^3 C_{4}^2+2016 C_{31}^2 C_{22} C_{4}^2+720 C_{31}^2 C_{4}^3+1408 C_{31} C_{22} C_{4}^3+512 C_{31} C_{4}^4+384 C_{22} C_{4}^4+128 C_{4}^5)^{12} (27648 C_{1111}^2 C_{22}^3+48384 C_{1111}^2 C_{22}^2 C_{4}+28224 C_{1111}^2 C_{22} C_{4}^2+5488 C_{1111}^2 C_{4}^3-13824 C_{1111} C_{211}^2 C_{22}^2-16128 C_{1111} C_{211}^2 C_{22} C_{4}-4704 C_{1111} C_{211}^2 C_{4}^2+5184 C_{1111} C_{211} C_{31}^2 C_{22}+3024 C_{1111} C_{211} C_{31}^2 C_{4}-20736 C_{1111} C_{211} C_{31} C_{22}^2-13824 C_{1111} C_{211} C_{31} C_{22} C_{4}-1008 C_{1111} C_{211} C_{31} C_{4}^2-2304 C_{1111} C_{211} C_{22}^2 C_{4}+2496 C_{1111} C_{211} C_{22} C_{4}^2+2240 C_{1111} C_{211} C_{4}^3-324 C_{1111} C_{31}^4+3888 C_{1111} C_{31}^3 C_{22}+972 C_{1111} C_{31}^3 C_{4}-5184 C_{1111} C_{31}^2 C_{22}^2+2160 C_{1111} C_{31}^2 C_{22} C_{4}+1080 C_{1111} C_{31}^2 C_{4}^2-6912 C_{1111} C_{31} C_{22}^3-8640 C_{1111} C_{31} C_{22}^2 C_{4}+1728 C_{1111} C_{31} C_{22} C_{4}^2+1280 C_{1111} C_{31} C_{4}^3-2304 C_{1111} C_{22}^3 C_{4}-1536 C_{1111} C_{22}^2 C_{4}^2+992 C_{1111} C_{22} C_{4}^3+320 C_{1111} C_{4}^4+1728 C_{211}^4 C_{22}+1008 C_{211}^4 C_{4}-144 C_{211}^3 C_{31}^2+5184 C_{211}^3 C_{31} C_{22}+2736 C_{211}^3 C_{31} C_{4}+576 C_{211}^3 C_{22} C_{4}+192 C_{211}^3 C_{4}^2-324 C_{211}^2 C_{31}^3+6912 C_{211}^2 C_{31}^2 C_{22}+3348 C_{211}^2 C_{31}^2 C_{4}+1728 C_{211}^2 C_{31} C_{22}^2+5472 C_{211}^2 C_{31} C_{22} C_{4}+2208 C_{211}^2 C_{31} C_{4}^2+576 C_{211}^2 C_{22}^2 C_{4}+1824 C_{211}^2 C_{22} C_{4}^2+832 C_{211}^2 C_{4}^3-324 C_{211} C_{31}^4+3888 C_{211} C_{31}^3 C_{22}+1404 C_{211} C_{31}^3 C_{4}+3456 C_{211} C_{31}^2 C_{22}^2+6912 C_{211} C_{31}^2 C_{22} C_{4}+2016 C_{211} C_{31}^2 C_{4}^2+2880 C_{211} C_{31} C_{22}^2 C_{4}+4368 C_{211} C_{31} C_{22} C_{4}^2+1184 C_{211} C_{31} C_{4}^3+960 C_{211} C_{22}^2 C_{4}^2+1088 C_{211} C_{22} C_{4}^3+224 C_{211} C_{4}^4-81 C_{31}^5+648 C_{31}^4 C_{22}+162 C_{31}^4 C_{4}+1296 C_{31}^3 C_{22}^2+1836 C_{31}^3 C_{22} C_{4}+432 C_{31}^3 C_{4}^2+576 C_{31}^2 C_{22}^3+2160 C_{31}^2 C_{22}^2 C_{4}+1800 C_{31}^2 C_{22} C_{4}^2+364 C_{31}^2 C_{4}^3+576 C_{31} C_{22}^3 C_{4}+1344 C_{31} C_{22}^2 C_{4}^2+800 C_{31} C_{22} C_{4}^3+128 C_{31} C_{4}^4+192 C_{22}^3 C_{4}^2+304 C_{22}^2 C_{4}^3+128 C_{22} C_{4}^4+16 C_{4}^5)^4$
\smallskip\\
\end{center}
An interesting particular case is

\begin{align}
S(x_1, x_2, x_3, x_4) = u (x_1 x_2 + x_1 x_3 + x_1 x_4 + x_2 x_3 + x_2 x_4 + x_3 x_4)^2 + v x_1x_2x_3x_4
\end{align}
\smallskip\\
In the context of Finslerian gravity, this quartic form is said to interpolate between the Minkovsky spacetime (described for certain reasons which we do not discuss here by the form $x_1 x_2 + x_1 x_3 + x_1 x_4 + x_2 x_3 + x_2 x_4 + x_3 x_4$) and the Berwald-Moor spacetime (described by $x_1x_2x_3x_4$). We can easily calculate its discriminant, using our formula (\ref{4Main}), and it turns out to be identically zero:

\begin{align*}
D_{4|4} \Big( u (x_1 x_2 + x_1 x_3 + x_1 x_4 + x_2 x_3 + x_2 x_4 + x_3 x_4)^2 + v x_1x_2x_3x_4 \Big) =
\end{align*}

\begin{align}
= {\cal D}_{4|4}\left( -\dfrac{v}{4}, \dfrac{v}{3}, \dfrac{u}{4} + \dfrac{v}{8}, - \dfrac{u}{2} - \dfrac{v}{4}, \dfrac{u}{4} + \dfrac{v}{24}\right) = 0
\end{align}
\smallskip\\
It is also instructive to consider the low-dimensional case $n = 2$. Formula (\ref{4Main}) gives

\begin{center}
${\cal D}_{2|4}\big(C_{4},  C_{31},  C_{22},  C_{211},  C_{1111}\big) = 2^{-6} C_{4}^2 (8 C_{1111}+4 C_{211}+2 C_{31}+2 C_{22}+C_{4}) (2 C_{22}+C_{4}) (-16 C_{1111} C_{22}-8 C_{1111} C_{4}+4 C_{211}^2+12 C_{211} C_{31}+8 C_{211} C_{4}+9 C_{31}^2+8 C_{31} C_{22}+16 C_{31} C_{4}+8 C_{22} C_{4}+8 C_{4}^2)^2$
\end{center}
Since this resultant is taken in $n = 2$ variables, one can alternatively perform a conventional Sylvester matrix calculation. The results precisely agree with each other. Just like in the previous section, the formula turns out to be valid for $n = 2$. This supports the hypothesis, that the main formula is actually valid for $n < r$ as well, despite it was derived in assumption $n \geq r$.

\section{Antisymmetric polynomials}

Apart from the completely symmetric polynomials, there are many other types of polynomials with symmetries and one can hope that their discriminant $D_{n|r}(S)$ is also simple. Following this line of thinking, we consider here the opposite case to symmetric polynomials -- completely antisymmetric polynomials

\begin{align}
S(x_1, \ldots, x_i, \ldots, x_j, \ldots, x_n) = - S(x_1, \ldots, x_j, \ldots, x_i, \ldots, x_n)
\end{align}
\smallskip\\
This case turns out to be even simpler, than the previously considered case of completely symmetric polynomials. Namely, let us show that for $n > 2$ discriminant of an antisymmetric polynomial vanishes. This is easy to do, since any completely antisymmetric polynomial $S(x_1, \ldots, x_n)$ can be written as

\begin{align}
S(x_1, \ldots, x_n) = \Delta(x_1, \ldots, x_n) {\widetilde S}(x_1, \ldots, x_n)
\end{align}
\smallskip\\
where ${\widetilde S}(x_1, \ldots, x_n)$ is some completely symmetric polynomial, and

\begin{align}
\Delta(x_1, \ldots, x_n) = \prod\limits_{i < j} (x_i - x_j)
\end{align}
\smallskip\\
is the Vandermonde determinant. Accordingly, the system of derivatives has a form

\begin{align}
\dfrac{\partial S}{\partial x_i} = \dfrac{\partial \Delta}{\partial x_i} {\widetilde S} + \Delta \dfrac{\partial {\widetilde S}}{\partial x_i}
\end{align}
\smallskip\\
and for $n > 2$ has plenty of critical points: for example, $(x_1, \ldots, x_n) = (1, \ldots, 1)$ is trivially a critical point, since for $n > 2$ not only $\Delta(1, \ldots, 1) = 0$ but also $\partial_i \Delta(1, \ldots, 1) = 0$. For $n = 2$, however, the situation is less trivial. The antisymmetric polynomial of degree $r$ in $n = 2$ variables can be written as

\begin{align}
S(x_1, x_2) = (x_1 - x_2) {\widetilde S}(x_1, x_2)
\end{align}
\smallskip\\
where ${\widetilde S}(x_1, x_2)$ is a symmetric polynomial of degree $(r - 1)$. The system of derivatives has a form

\begin{align}
\dfrac{\partial S}{\partial x_1} = {\widetilde S} + (x_1 - x_2) \dfrac{\partial {\widetilde S}}{\partial x_1} \\
\nonumber \\
\dfrac{\partial S}{\partial x_2} = - {\widetilde S} + (x_1 - x_2) \dfrac{\partial {\widetilde S}}{\partial x_2}
\end{align}
\smallskip\\
which can be, using the Euler identity $x_1 \dfrac{\partial S}{\partial x_1} + x_2 \dfrac{\partial S}{\partial x_2} = r S$ rewritten as

\begin{align}
\dfrac{\partial S}{\partial x_1} = (x_1 + \dfrac{x_1}{r} - x_2) \dfrac{\partial {\widetilde S}}{\partial x_1} + \dfrac{x_2}{r} \dfrac{\partial {\widetilde S}}{\partial x_1} \\
\nonumber \\
\dfrac{\partial S}{\partial x_2} = - \dfrac{x_1}{r} \dfrac{\partial {\widetilde S}}{\partial x_1} + (x_1 - x_2 - \dfrac{x_2}{r}) \dfrac{\partial {\widetilde S}}{\partial x_2}
\end{align}
\smallskip\\
One can see, that $S$ can be degenerate only in two cases: either if ${\widetilde S}$ is degenerate, or if

\begin{align}
\det\limits_{2 \times 2} \left( \begin{array}{cc} x_1 + \dfrac{x_1}{r} - x_2 & \dfrac{x_2}{r} \\ \\ - \dfrac{x_1}{r} &  x_1 - x_2 - \dfrac{x_2}{r} \end{array} \right) = \left( 1 + \dfrac{1}{r} \right) (x_1 - x_2)^2 = 0
\end{align}
\smallskip\\
In the second case $x_1 = X$, $x_2 = X$ and $\partial S / \partial x_1 = - \partial S / \partial x_2 = {\widetilde S}(1,1) X^{r-1}$. Consequently, $S$ can be degenerate either if ${\widetilde S}$ is degenerate, or if ${\widetilde S}(1,1) = 0$. In terms of discriminants this statement takes form

\begin{align}
D_{2|r}(S) = {\widetilde S}(1,1)^2 \cdot D_{2|r-1}({\widetilde S})
\end{align}
\smallskip\\
One can see that discriminants of antisymmetric polynomials are vanishing for $n > 2$, and for $n = 2$ they are expressed through discriminants of symmetric polynomials. The case of $n = 2$ is, however, of little interest, since conventional Sylvester formula is quite enough for any practical purpose.

\section{Summary}

Any homogeneous polynomial $S(x_1, \ldots, x_n)$ of degree $r$ in $n$ variables posesses a discriminant $D_{n|r}(S)$, which is typically a complicated polynomial in the coefficients of $S$. For generic $S$, writing a closed expression for $D_{n|r}(S)$ is almost hopeless: already the discriminant $D_{3|3}(S)$ contains 2040 monomials (and takes several pages). We have shown, that in particular case of symmetric polynomials

\begin{align}
S(x_1, \ldots, x_i, \ldots, x_j, \ldots, x_n) = S(x_1, \ldots, x_j, \ldots, x_i, \ldots, x_n)
\end{align}
\smallskip\\
discriminant greatly simplifies and becomes a product of many smaller factors:

\begin{align*}
{\cal D}_{n|r} \left( S \right) \ \sim \ \prod\limits_{M_1 + \ldots + M_{r-1} = n} \Big( d_{M}\left( S \right) \Big)^{ \dfrac{\#_M!}{(r-1)!} \dfrac{(M_1 + \ldots + M_{r-1})!}{M_1! \ldots M_{r-1}!}  }
\end{align*}
\smallskip\\
Because of this factorisation, discriminants of symmetric polynomials do not take hundreds and hundreds of pages: much shorter closed expressions for them can be written out. Each particular factor $d_{M}\left( S \right)$ can be computed as a resultant in no more than in $(r - 1)$ variables. This allows to calculate discriminants of symmetric polynomials even for $n >> r$. It may be even possible to study different kinds of large-$n$ limits, corresponding to infinite-dimensional or functional discriminants. Interesting directions of generalisation are to other types of symmetries and to non-homogeneous polynomials.

\section{Acknowledgements}

We are indebted to A.Morozov for illuminating discussions and kind support. Our work is partly supported by Russian Federal Nuclear Energy Agency and the Russian President's Grant of Support for the Scientific Schools NSh-3035.2008.2, by RFBR grant 07-02-00645, by the joint grants 09-01-92440-CE, 09-02-91005-ANF and 09-02-93105-CNRS. The work of Sh.Shakirov is also supported in part by the Moebius Contest Foundation for Young Scientists.

\end{document}